\documentclass[11pt]{amsart}
\usepackage{amsmath,amssymb}
\usepackage{graphicx}
\usepackage{xcolor}
\usepackage{xkeyval}
\usepackage{calc}
\usepackage{ifthen}
\usepackage{algorithm}
\usepackage{algpseudocode}
\usepackage{mathrsfs}
\usepackage{fullpage}
\usepackage[pdftex,colorlinks,backref=page,citecolor=blue]{hyperref}
\usepackage{enumerate}
\usepackage{pgf,tikz}
\usepackage{todonotes}
\usepackage{appendix}
\usetikzlibrary{arrows}
\usetikzlibrary{positioning}
\usepackage{array}

\def\COMMENT#1{}
\def\TASK#1{}

\newdimen\margin   
\def\textno#1&#2\par{%
    \margin=\hsize
    \advance\margin by -4\parindent
           \setbox1=\hbox{\sl#1}%
    \ifdim\wd1 < \margin
       $$\box1\eqno#2$$%
    \else
       \bigbreak
       \hbox to \hsize{\indent$\vcenter{\advance\hsize by -3\parindent
       \sl\noindent#1}\hfil#2$}%
       \bigbreak
    \fi}

\makeatletter \providecommand\@dotsep{5} \def\listtodoname{List of Todos} \def\listoftodos{\@starttoc{tdo}\listtodoname} \makeatother

\newtheorem{thm}{Theorem}[section]
\newtheorem{define}[thm]{Definition}

\newtheorem{lem}[thm]{Lemma}

\newtheorem{cor}[thm]{Corollary}

\newtheorem{prop}[thm]{Proposition}

\newtheorem*{thm*}{Theorem}
\newtheorem*{define*}{Definition}
\newtheorem*{examp*}{Example}
\newtheorem*{lem*}{Lemma}
\newtheorem*{claim*}{Claim}
\newtheorem*{fact*}{Fact}
\newtheorem*{col*}{Corollary}
\newtheorem*{conj*}{Conjecture}

\title{A degree sequence strengthening of the vertex degree threshold for a perfect matching in $3$-uniform hypergraphs}

\author{Candida Bowtell \and Joseph Hyde}\thanks{CB: Mathematics Institute, University of Warwick, United Kingdom, {\tt Candy.Bowtell@warwick.ac.uk}, research supported in part by the European Research Council (ERC) under the European Union Horizon 2020 research and innovation programme (grant agreement No. 947978).\\
\indent JH: Mathematics and Statistics, University of Victoria, Victoria, B.C., Canada, {\tt josephhyde@uvic.ca}, research supported
by the UK Research and Innovation Future Leaders Fellowship MR/S016325/1.\\
}

\begin{document}
\maketitle

\begin{abstract}
The study of asymptotic minimum degree thresholds that force matchings and tilings in hypergraphs is a lively area of research in combinatorics. A key breakthrough in this area was a result of H\`{a}n, Person and Schacht \cite{hps} who proved that the asymptotic minimum vertex degree threshold for a perfect matching in an $n$-vertex $3$-graph is $\left(\frac{5}{9}+o(1)\right)\binom{n}{2}$. In this paper we improve on this result, giving a family of degree sequence results, all of which imply the result of H\`{a}n, Person and Schacht, and additionally allow one third of the vertices to have degree $\frac{1}{9}\binom{n}{2}$ below this threshold. Furthermore, we show that this result is, in some sense, tight.
\end{abstract}

\section{Introduction}

Determining whether a $k$-uniform hypergraph (or $k$-graph) contains a perfect matching (a collection of disjoint edges which cover the vertex set) is a key question in combinatorics. Whilst Tutte's Theorem \cite{tutte} gives a complete characterisation for when a graph $G$ contains a perfect matching, no such `nice' characterisation is expected to be found for $k$-graphs in general. In particular, determining whether a $k$-graph contains a perfect matching is one of Karp's original 21 NP-complete problems \cite{karp}. As such, much work has been done to consider sufficient conditions for a hypergraph to contain a perfect matching. A key direction for this has been to consider minimum degree conditions, also known as `Dirac-type' conditions, which follow the form of Dirac's Theorem \cite{d} from 1952; every graph $G$ on $n \geq 3$ vertices with minimum degree at least $n/2$ contains a Hamilton cycle (a cycle covering all vertices in $G$), and if $n$ is even, contains a perfect matching (found by taking every other edge in a Hamilton cycle). In hypergraphs, the notion of minimum degree extends in various ways. In particular, for a hypergraph $H$, we define the \emph{degree of a set $T \subseteq V(H)$, $\deg(T)$,} to be the number of edges in $H$ containing $T$. We then define the \emph{minimum $t$-degree, $\delta_t(H)$,} to be  $\delta_t(H):=\min\{\deg(T): T \subseteq V(H), |T|=t\}$. In a $k$-graph, $H$, we also refer to $\delta_1(H)$ as the \emph{minimum vertex degree of $H$}, and to $\delta_{k-1}(H)$ as the \emph{minimum co-degree of $H$}. In the last fifteen years much progress has been made in finding asymptotic and exact minimum $t$-degree conditions that force perfect matchings in $k$-graphs for various $k$ and $t$. Let $m_t(k,n)$ denote the smallest integer $m$ such that every $k$-graph on $n$ vertices with minimum $t$-degree at least $m$ contains a perfect matching (given, of course, also that $n \in k\mathbb{Z}$). We refer to  $m_t(k,n)$ as the minimum $t$-degree threshold for a $k$-graph on $n$ vertices to contain a perfect matching. For the purposes of this paper, we are interested in minimum vertex degree thresholds, and particularly, in $m_1(3,n)$. For more on Dirac-type problems in general, see e.g. \cite{dirac_survey1} and \cite{dirac_survey2}. In a significant breakthrough, H\`{a}n, Person and Schacht \cite{hps} determined the asymptotic minimum vertex degree threshold for a perfect matching in a $3$-graph:

\begin{thm}[H\`{a}n, Person and Schacht, \cite{hps}] \label{thm:hps_main}
For all $\gamma>0$ there exists an $n_0 = n_0(\gamma) \in \mathbb{N}$ such that for all $n \geq n_0$ with $n \in 3\mathbb{Z}$ the following holds. Let $H$ be a $3$-graph on $n$ vertices with
$$\delta_1(H) \geq \left(\frac{5}{9}+\gamma\right)\binom{n}{2}.$$
Then $H$ contains a perfect matching.
\end{thm}

This was subsequently improved to an exact result by K\"{u}hn, Osthus and Treglown \cite{kot_3_match}, and independently by Khan \cite{khan}; that is, $m_1(3,n)=\binom{n-1}{2}-\binom{2n/3}{2}+1$ for sufficiently large $n$. We can see this result is tight by examining the following extremal example: Let $H$ be a $3$-graph on $n$ vertices and divide $V(H)$ into two parts, $A$ and $B$, with $|A|=\frac{n}{3} -1$ and $|B|=\frac{2n}{3}+1$. Let $E(H)$ contain all edges with at least one vertex in $A$. Observe that $H$ does not contain a perfect matching, since every edge in a matching will have at least one vertex in $A$, and a perfect matching has size $n/3$, but $A$ has only $n/3-1$ vertices. This extremal example can be generalised from $3$-graphs to $k$-graphs, and is known as the \emph{space barrier}, a term coined by Keevash and Mycroft \cite{keev-mycroft}. 

Whilst asymptotic and exact results for $m_t(k,n)$ are best possible in the sense that one cannot lower the minimum $t$-degree threshold and still ensure the existence of a perfect matching, we can also consider `stronger' degree conditions, by seeing whether allowing a certain proportion of $t$-sets to go a certain distance below the minimum $t$-degree threshold, we are still able to guarantee a perfect matching. This idea is formalised by the notion of a \emph{degree sequence} of a graph. In particular, we say that a graph $G$ on $n$ vertices has degree sequence $d_1 \leq d_2 \leq \ldots \leq d_n$ if there exists an ordering $(v_1, v_2, \ldots, v_n)$ of the vertices of $G$ such that $d(v_i) = d_i$ for all $i \in [n]$. It is natural to ask for which degree sequences of $G$ we are guaranteed a perfect matching. In general, it is hard to characterise all such degree sequences, but there are notable results of P\'{o}sa \cite{posa} and Chv\'{a}tal \cite{chv} which show two different degree sequence improvements of Dirac's Theorem. P\'{o}sa \cite{posa} proved that if $G$ is a graph on $n \geq 3$ vertices with degree sequence $d_1 \leq \ldots \leq d_n$ satisfying
$d_i \geq i+1$ for all $i < (n-1)/2$ and, when $n$ is odd $d_{\lceil n/2 \rceil} \geq \lceil n/2 \rceil$, then $G$ contains a Hamilton cycle. Chv\'{a}tal \cite{chv} went further and demonstrated that if $G$ has degree sequence $d_1 \leq d_2 \leq \ldots \leq d_n$ such that either $d_i \geq i+1$ or $d_{n-i} \geq n-i$ for all $i \leq n/2$ then $G$ still contains a Hamilton cycle. Moreover, for every sequence not satisfying this condition, there is a graph with pointwise at least as
large degree sequence not containing a Hamilton cycle.
In this paper we will be concerned with so-called \emph{P\'{o}sa-type} degree sequence conditions, that is, degree sequence conditions which, informally, have a starting point ($d_1$) below some (known) minimum degree threshold and have a part of the degree sequence condition that steadily increases.

There have been a number of recent examples of P\'{o}sa-type degree sequence results in graphs. Asymptotically answering a conjecture of Balogh, Kostochka and Treglown \cite{bkt}, Treglown \cite{Tregs} proved a P\'{o}sa-type degree sequence version of Hajnal and Szemer\'{e}di's \cite{hs} perfect $K_r$-tiling result (as well as a degree sequence strengthening of Alon and Yuster's \cite{ay} perfect $H$-tiling result for general graphs $H$). 
Using ideas from \cite{Tregs}, Hyde, Liu and Treglown \cite{hlt} proved a P\'{o}sa-type degree sequence strengthening of Koml\'{o}s' \cite{Komlos} almost-perfect tiling theorem which was then utilised by Hyde and Treglown \cite{ht} to give a P\'{o}sa-type degree sequence version of K\"{u}hn and Osthus' \cite{kuhn} perfect tiling theorem. 
See \cite{kt, ps, st} for further examples of degree sequence results.

Whilst considerable progress has been made with respect to degree sequence results in graphs, not as much headway has been made for hypergraphs. Very recently, Sch\"{u}lke \cite{schulke} proved a P\'{o}sa-type degree sequence result related to finding tight Hamilton cycles in $3$-graphs. We say that a 3-graph $H$ on $n$ vertices contains a tight Hamilton cycle if there exists an ordering  $(v_1, v_2, \ldots, v_n)$ of the vertices of $H$ such that $\{v_iv_{i+1}v_{i+2}, \ i \in [n-2]\} \cup v_{n-1}v_{n}v_{1} \cup v_{n}v_{1}v_{2} \subseteq E(H)$. Furthermore, for $n \in \mathbb{N}$ and a 3-graph $H = ([n], E)$ we define $d(i,j)$ to be the number of edges of $H$ containing both vertex $i$ and vertex $j$. Generalising a result of R\"{o}dl, Ruci\'{n}ski and Szemer\'{e}di \cite{rrs} on the asymptotic minimum co-degree threshold for a tight Hamilton cycle in a 3-graph, Sch\"{u}lke \cite{schulke} proved the following:

\begin{thm}\label{thm:schulke}
For all $\gamma > 0$ there exists an $n_0=n_0(\gamma) \in \mathbb{N}$ such that for all $n \in \mathbb{N}$ with $n \geq n_0$ the following holds. If $H = ([n], E)$ is a 3-graph with $d(i,j) \geq \min\left\{i,j, \frac{n}{2}\right\} + \gamma n \ \mbox{for all} \ \{i,j\} \in \binom{[n]}{2},$ then $H$ contains a (tight) Hamilton cycle. 
\end{thm}

Theorem~\ref{thm:schulke} can be seen as an analogue of P\'{o}sa's theorem for $3$-graphs. The proof follows the strategy taken in \cite{rrs}. (Note that Theorem~\ref{thm:schulke} yields a perfect matching in $H$ whenever $n \in 3\mathbb{Z}$ by taking every third edge in a Hamilton cycle.) We believe that Theorem~\ref{thm:schulke} is the first sufficient degree sequence condition for the existence of some spanning structure in a hypergraph. In particular, as far as we are aware, no work has been done to provide degree sequence improvements to minimum vertex degree thresholds for structures in $k$-graphs. Note that for both the graph case and the co-degree case, the largest possible degree in a $k$-graph is $n-k+1$, and so a degree sequence result typically has a gap of $\Theta(n)$ between the smallest and largest degrees in the degree sequence condition. 
However, a substantial difference for $t$-degree conditions in $k$-graphs, where $t<k-1$, is that, to make a significant improvement on the minimum degree threshold, we wish to have the minimum $t$-degree in the degree sequence starting a constant proportion lower than the minimum $t$-degree threshold, which typically means increasing the degree by $\Theta(n^{k-t})$ (where $k-t \geq 2$). Our main result, a collection of P\'{o}sa-type degree sequence strengthenings of Theorem~\ref{thm:hps_main}, is the following:

\begin{thm}\label{thm:main}
For all $\gamma > 0$ there exists $n_0 = n_0(\gamma) \in \mathbb{N}$ such that for every $n \geq n_0$ with $n \in 3\mathbb{N}$ and $q \in [(1-\sqrt{\frac{2}{3}})n]$, the following holds. Suppose $H$ is a 3-graph on $n$ vertices with degree sequence $d_1 \leq \ldots \leq d_n$ such that 
$$d_i \geq \begin{cases}
    \left(\frac{1}{3} + \gamma \right) \binom{n}{2} + iq  & \mbox{if}\  1 \leq i \leq q,\\
    \left(\frac{4}{9} + \gamma \right)\binom{n}{2}  & \mbox{if}\  q < i \leq \frac{n}{3}, \\
   \left(\frac{5}{9} + \gamma \right)\binom{n}{2} &  \mbox{if}\  \frac{n}{3} < i. 
    \end{cases}$$ 
Then $H$ contains a perfect matching.
\end{thm}

Note that our result contains many (growing with $n$) degree sequence conditions which do not imply one another, all of which imply Theorem~\ref{thm:hps_main}, and improve on Theorem~\ref{thm:hps_main} by having a third of the vertices substantially below the minimum degree threshold in Theorem~\ref{thm:hps_main}. Theorem~\ref{thm:main} is tight in the sense that it is not possible to have more than a third of the vertices below the $4/9$ barrier by any $\omega(n)$ amount,\footnote{That is, there exists a constant $c>0$ such that it is not possible to have more than a third of the vertices with degree less than $(1-\frac{c}{n})\frac{4}{9}\binom{n}{2}$. In general, from now on, when we refer to the relation of vertices and their degree to `the $x$ barrier' we mean their relation to $(1+o(1))x\binom{n}{2}$ where sometimes, as here, we are more precise in the $o(1)$ term.} as seen by the following space and parity examples.

\subsection{Extremal example 1: space barrier}
Let $H$ be a $3$-graph with vertex set $V(H)=A~\dot\cup~B$, where $|A|=\frac{n}{3}+1$, $|B| = \frac{2n}{3} - 1$ and $E(H)$ consists of all edges containing at most one vertex from $A$. Then for each $v \in A$, 
$$\deg(v)=\binom{\frac{2n}{3}-1}{2}=\frac{2}{9}n^2-n+1=\frac{4}{9}\binom{n}{2}-\frac{7n}{9}+1,$$
and for each $v \in B$, 
$$\deg(v)=\binom{n-1}{2}-\binom{\frac{n}{3}+1}{2} \sim \frac{8}{9}\binom{n}{2}.$$
$H$ has no perfect matching, since each of the $n/3+1$ vertices in $A$ uses two vertices in $B$ to form an edge, and so to cover $A$ we need at least $2(n/3+1)$ vertices in $B$. Hence we cannot cover $A$. This implies that we cannot have more than $\frac{n}{3}$ vertices with degree $\omega(n)$ below the $4/9$ barrier.

\subsection{Extremal example 2: parity barrier}
Let $H$ be a $3$-graph with vertex set $V(H)=A \dot\cup B$, such that $|A|=\frac{n}{3}$ is odd (and $|B| = \frac{2n}{3}$). Let $E(H)$ consist of all edges with an even number of vertices in $A$. Then for $v \in A$,
$$\deg(v)=\left(\frac{n}{3}-1\right)\frac{2n}{3}= \frac{2n^2}{9}-\frac{2n}{3}=\frac{4}{9}\binom{n}{2}-\frac{4n}{9}.$$
Clearly $H$ has no perfect matching, because every edge using a vertex from $A$ has to use exactly two vertices from $A$, so since $|A|$ is odd, it is not possible to cover $A$ by disjoint edges. Note that in this example we have that every vertex $v \in B$ has $\deg(v)=m_1(3,n)-1$. Considering this, it is perhaps not so surprising that this does not have a perfect matching. However, if we instead take $|A|=\frac{n}{3}-1$ (and still require that $|A|$ is odd), we get that for vertices $v \in A$, $\deg(v)=\frac{4}{9}\binom{n}{2}-\frac{7n}{9}-2$ and that for vertices $v \in B$, $\deg(v) \geq m_1(3,n)$ and again there is no perfect matching. Off the back of these two cases combined, we wonder whether it would be possible to prove the following statement: every $3$-graph $H$ on $n \in 3\mathbb{N}$ vertices such that $n$ is sufficiently large and $2n/3$ vertices in $H$ have degree at least $m_1(3,n)$ whilst the remaining $n/3$ have degree at least $\frac{2n^2}{9}- \frac{2n}{3} +1$ has a perfect matching. If this were true, the above parity barrier would be tight since it shows that having each vertex degree only one lower results in no perfect matching. (This is tight in the same sense that $m_1(3,n)$ is a tight threshold - the extremal example there shows that we do not get a perfect matching if $\frac{2n}{3}+1$ vertices have degree only one below that threshold.)
\smallskip

It seems difficult to find extremal examples to suggest that the sequences in Theorem~\ref{thm:main} are exactly optimal. On the other hand, as we discuss in Section~\ref{sec:comments}, new ideas would be needed to potentially improve our result. We hope that our result will lead to further development of the area of hypergraph degree sequences, and further understanding of the variety of degree sequence improvements possible for known minimum degree threshold results in hypergraphs.

\subsection{Paper organisation}

Our proof of Theorem~\ref{thm:main} is split into two parts: an almost-perfect matching lemma and an absorbing lemma. The proof of the former employs ideas from \cite{hlt} and \cite{Tregs}, in particular with regards to the `swapping' arguments employed in \cite{hlt}, and also relies on inferences one can make from the proof of \cite[Theorem 4.4]{hps} (see the proof of Lemma~\ref{lem:largermatch}). The proof of the latter borrows from the proof of \cite[Lemma 2.4]{hps}, with some new ideas introduced to accommodate vertices with small degree.

The rest of the paper is laid out as follows: in Section~\ref{sec:apm}, we discuss the statement of our almost-perfect matching lemma, and its relation to \cite[Theorem 4.4]{hps}. In Section~\ref{sec:details} we introduce the notion of swapping pairs, and present the key details of the proof to obtain our almost-perfect matching (Theorem~\ref{thm:almostmain}), and in Section~\ref{sec:absorbing} we discuss our modified absorbing argument and complete the proof of Theorem~\ref{thm:main}. Finally, in Section~\ref{sec:comments}, we discuss directions for future development.

\subsection{Preliminary definitions and notation}

We write $[n]:=\{1, 2, \ldots, n\}$. For $l \in \mathbb{N}$ and a collection of sets $A$, we let $\binom{A}{l}:=\{S \subseteq A: |S| = l\}$, that is, $\binom{A}{l}$ contains the unordered $l$-sets of elements from $A$, not allowing repeats. We write $A^{(l)}$ to denote the collection of unordered $l$-sets of elements from $A$, where repeats {\bf are} allowed.

We define a \emph{$3$-graph} $H$ to be a set of vertices $V(H)$ together with an edge set $E(H)$ consisting of $3$-sets of vertices from $V(H)$. Let $X \subseteq V(H)$. Then $H[X]$ is the \textit{subhypergraph of $H$ induced by $X$} and has vertex set $X$ and edge set $E(H[X]) := \{xyz \in E(H): x,y,z \in X\}$. We also define $H\setminus X = H[V(H)\setminus X]$. For a set $M \subseteq \binom{V(H)}{l}$, we write $V(M) := \bigcup_{m \in M}m$.

For a 3-graph $H$ on $n$ vertices with degree sequence $d_1 \leq \ldots \leq d_n$ and some given $\gamma \in \mathbb{R}$, we partition the vertex set $V(H)$ into three families according to their position in the degree sequence.
Let $V_{5/9}(H,\gamma):=\{v \in V(H): d(v) \geq (\frac{5}{9}+\gamma) \binom{n}{2}\}$, $V_{4/9}(H,\gamma):= \{v \in V(H): d(v) \geq (\frac{4}{9}+\gamma) \binom{n}{2}\} \setminus V_{5/9}(H,\gamma)$, and  $V_{3/9}(H,\gamma):=\{v \in V(H): d(v) < (\frac{4}{9}+\gamma) \binom{n}{2}\}$. We write $V_{5/9}$, $V_{4/9}$ and $V_{3/9}$, respectively, when $n$ and $\gamma$ are clear from context. We refer to the vertices in these sets as $5/9$th, $4/9$th and $3/9$th vertices respectively, and also as big, medium and small vertices, respectively. Furthermore, we say that any vertex in $V_{4/9} \cup V_{5/9}$ is \emph{not-small}.

\section{The almost-perfect matching} \label{sec:apm}

For a matching $M$ in a 3-graph $H$, let $L(M):= V(H) \setminus V(M)$ be the \emph{leave of $M$}. Most of our work concerns the proof of the following theorem:

\begin{thm}\label{thm:almostmain}
Let $\gamma > 0$. There exists $n'' = n''(\gamma) \in \mathbb{N}$ such that for every $n \geq n''$ and $q \in \left[\frac{n}{3\sqrt{2}}\right]$, the following holds. Let $H$ be a $3$-graph on $n$ vertices with degree sequence $d_1 \leq \ldots \leq d_n$ such that
$$d_i \geq \begin{cases}
    \left(\frac{1}{3}+4\gamma \right)\binom{n}{2} + iq & \mbox{if}\  1 \leq i \leq q,\\
    \left(\frac{4}{9}+4\gamma \right)\binom{n}{2} & \mbox{if}\  q< i \leq \frac{n}{3},\\
    \left(\frac{5}{9}+4\gamma \right)\binom{n}{2} &  \mbox{if}\  \frac{n}{3} < i.
\end{cases}$$ 
Then $H$ contains a matching $M$ of size $\lfloor \frac{n-\gamma n}{3} \rfloor$ where $|L(M) \cap V_{5/9}| \geq \frac{2}{3}|L(M)|$.
\end{thm}

Note that we have $q \leq \frac{n}{3\sqrt{2}}$, since $\left(\frac{1}{3}+4\gamma\right)\binom{n}{2}+q^2 > \left(\frac{4}{9}+4\gamma\right)\binom{n}{2}$ when $q > \frac{n}{3\sqrt{2}}$, and so if we take a larger value for $q$ we obtain a pointwise larger degree sequence than that in Theorem~\ref{thm:almostmain}.

In order to prove this theorem, we both have to show that there is a matching $M$ of the required size, and that $L(M)$ contains sufficiently many $5/9$th vertices.
The key strategy for the proof is the use of a swapping mechanism to show that if we do not have enough $5/9$th vertices in the leave of the current matching, we can find a matching of at least the same size, which increases the number of $5/9$th vertices in the leave. In particular, we first show that, given a largest matching which has size at most $\lfloor \frac{n-\gamma n}{3} \rfloor$, we can find a matching of the same size with a constant proportion of big vertices in the leave (see Lemma~\ref{lem:betterdegree}). Then we can infer from the proof of \cite[Theorem 4.4]{hps} that a larger matching exists, contradicting that our first choice of matching was largest. Hence there exists a matching of size $\lfloor \frac{n-\gamma n}{3} \rfloor$ (see Lemma~\ref{lem:largermatch}). Once we have such a matching, we again use the swapping arguments to show that we may now obtain a matching of the same size with the required proportion of big vertices in the leave (see proof of Theorem~\ref{thm:almostmain} on p.\pageref{proof:almostmainthm}).

The following definition is crucial in the proof of \cite[Theorem 4.4]{hps} and for our subsequent swapping arguments.

\begin{define}
\textnormal{Let $H$ be a 3-graph, let $U$ be a collection of disjoint $3$-sets of vertices in $V(H)$ and let $v \in V(H) \setminus V(U)$. We define the \emph{$3$-set link graph} $L_v(U)$, to be the (2-)graph on vertex set 
$$V(L_v(U)) := \underset{u \in U}{\bigcup} u,$$ 
and edge set given by any pair of vertices from distinct $3$-sets in $U$ that together with $v$ form an edge in $H$. i.e. $uw \in E(L_v(U))$ if and only if there exist $e, f \in U$ with $e \neq f$ such that $u \in e, w \in f$ and $vuw \in E(H)$. By slight abuse of notation, for $e,f \in \binom{V(H)}{3}$, we sometimes write $L_v(e,f)$ in place of $L_v(\{e,f\})$.}
\end{define}

\begin{lem}\label{lem:largermatch}
Suppose that $H$ is as in Theorem~\ref{thm:almostmain} and that $M$ is a matching in $H$ such that $|M| < \lfloor \frac{n-\gamma n}{3}\rfloor$ and $|L(M) \cap V_{5/9}| > \frac{2\gamma n}{75}$. Then there exists a matching $M^*$ in $H$ with $|M^*| \geq |M|+1$.
\end{lem}
\begin{proof}
Suppose for a contradiction that no such $M^*$ exists. Let $B(M) \subseteq L(M)$ be the set of $5/9$th vertices in $L(M)$. By assumption we have that $|B(M)|> \frac{2\gamma n}{75}$. Let $s:=\lfloor \frac{n-\gamma n}{3}\rfloor - |M|$. Then take $S \subseteq \binom{L(M)}{3}$ with $|S| = s$ and $|B(M) \setminus V(S)| > \frac{2\gamma n}{75}$ (i.e. take as many non-big vertices from $L(M)$ for $S$ before adding any big vertices). Let $N=M \cup S$ so that $|N|=\lfloor \frac{n-\gamma n}{3}\rfloor$. Then for $L(N):=V(H) \setminus V(N)$, and $B(N):=B(M) \cap L(N)$, we still have that $|B(N)| >  \frac{2\gamma n}{75}$.  Following the proof of \cite[Theorem 4.4]{hps},\footnote{Appendix~\ref{app1} includes a brief summary of the strategy used in the proof of \cite[Theorem 4.4]{hps}, and the key details we take from it.} we first note that for every $v \in B(N)$, we have 
$$|E(L_v(N))| \geq \deg_H(v) - 3|N| - |L(N)|(n-|L(N)|) - \binom{|L(N)|}{2}> \left(\frac{5}{9} + \gamma\right)\binom{n}{2}.$$
Their proof shows that either we may find a larger matching $|M^*|$, with $|M^*| \geq |M|+1$, or we have at most $\frac{2\gamma n}{75}$ vertices satisfying
$$|E(L_v(N))| > \left(\frac{5}{9} + \gamma\right)\binom{n}{2}$$
in $L(N)$. However, since $|B(N)|>\frac{2\gamma n}{75}$, we must be in the former case, and hence able to find the desired matching.
\end{proof}

The following lemma is the heart of our proof of Theorem~\ref{thm:almostmain} and will be proved in the next section.

\begin{lem} \label{lem:betterdegree}
Let $H$ be as in Theorem~\ref{thm:almostmain} and $M$ be a matching of $H$ such that $|M| \leq \lfloor \frac{n-\gamma n}{3}\rfloor$. Then there exists a matching $N$ such that $|N| \geq |M|$ and $|L(N) \cap V_{5/9}| > \frac{2\gamma n}{75}$.
\end{lem}

These lemmas imply the following key corollary:

\begin{cor} \label{cor:betterdegree}
For $H$ as in Theorem~\ref{thm:almostmain} there exists a matching $M^*$ in $H$ such that $|M^*|=\lfloor \frac{n-\gamma n}{3}\rfloor$ and $|L(M^*) \cap V_{5/9}| > \frac{2\gamma n}{75}$. 
\end{cor}




\section{Proof of Lemma~\ref{lem:betterdegree}} \label{sec:details}

In this section we introduce the key swapping lemmas which allow us to obtain the required matchings with sufficiently many 5/9th vertices in the leave.
Throughout this section for a given maximum matching $M'$ in our graph $H$, we define a \emph{phantom matching, $M$, of $M'$} in the following way: if $|M'| < \lfloor \frac{n-\gamma n}{3} \rfloor$ we define $M \supseteq M'$ such that $M$ consists of disjoint $3$-sets and $|M|=\lfloor \frac{n-\gamma n}{3} \rfloor$.
If $|M'| \geq \lfloor \frac{n-\gamma n}{3} \rfloor$ we define $M \subseteq M'$ such that $|M|=\lfloor \frac{n-\gamma n}{3} \rfloor$, so that our phantom matching always has size $\lfloor \frac{n-\gamma n}{3} \rfloor$.
We shall refer directly to a phantom matching $M$, meaning a collection of disjoint $3$-sets from $V(H)$ such that there exists a maximum matching $M'$ such that $M$ is a phantom matching of $M'$. Also, the $3$-sets in a phantom matching will sometimes be referred to as \emph{phantom edges}.

\subsection{Swapping pairs}

The swapping arguments we use require a detailed understanding of the different combinations of vertices that may reside in each phantom edge in a phantom matching. As such we have a substantial set of notation to deal with the different cases, which is explained here. Throughout this section, unless stated otherwise, we shall call the vertices of a $3$-set $e$ by $e_1, e_2, e_3$.

\begin{define}
\textnormal{Let $H$ be a 3-graph and $e, f \in \binom{V(H)}{3}$ be disjoint $3$-sets in $H$. Let $x,y \in V(H) \setminus (e \cup f)$. We say \emph{$\{e,f\}$ has an $\{x,y\}$-matching} if there exist vertices $e_1, e_2 \in e$ and $f_1, f_2 \in f$ with $e_1 \neq e_2$ and $f_1 \neq f_2$ such that $xe_1f_1, ye_2f_2 \in E(H)$. We call $\{xe_1f_1, ye_2f_2\}$ an \emph{$\{x,y\}$-matching for $\{e,f\}$}. Furthermore, in this subsection denote the vertices in $e$ and $f$ not present in this $\{x,y\}$-matching by $e_3$ and $f_3$ respectively.}
\end{define}

\begin{define} 
\textnormal{Let $H$ be a 3-graph on $n$ vertices with degree sequence $d_1 \leq \ldots \leq d_n$. Define a bijection $I_H: V(H) \to [n]$ such that $I_H(x) = i$ implies that $d_H(x) := d_i$. This natural bijection will be used several times throughout this section. Let $x$ and $y$ be $3/9$th vertices in $H$. We say that} 
$$y \ \mbox{is} \ 
\begin{cases}
x\mbox{-little} & \mbox{if} \ I_H(y) < I_H(x);\\
 x\mbox{-large} & \mbox{if} \ I_H(y) > I_H(x).
\end{cases}$$
\end{define}

\begin{define}
\textnormal{Let $H$ be a 3-graph and $M$ be a phantom matching in $H$. Let $x,y \in L(M)$ and $\{e,f\} \in \binom{M}{2}$. We say \emph{$\{x,y\}$ is a $5/9$th (or big) swapping pair for $\{e,f\}$} if there exists an $\{x,y\}$-matching for $\{e,f\}$ such that $\{e_3, f_3\} \cap V_{5/9} \neq \emptyset$. We say that \emph{$\{x,y\}$ is a $4/9$th (or not-small)  swapping pair for $\{e,f\}$} if there exists an $\{x,y\}$-matching for $\{e,f\}$ such that $\{e_3, f_3\} \cap (V_{4/9} \cup V_{5/9}) \neq \emptyset$. We say \emph{$\{x,y\}$ is a large swapping pair for $\{e,f\}$} if there exists an $\{x,y\}$-matching for $\{e,f\}$ such that both $e_3$ and $f_3$ are $x$-large and $y$-large. In general, if $\{x,y\}$ is a $5/9$th, $4/9$th or large swapping pair for $\{e,f\}$ then we say that \emph{$\{x,y\}$ is a swapping pair for $\{e,f\}$}. We call an $\{x,y\}$-matching for $\{e,f\}$ \emph{good} if it is a witness for $\{x,y\}$ being a swapping pair for $\{e,f\}$.}
\end{define}

For a phantom matching $M$ and a vertex $x \in L(M)$, we now describe different subsets of $M$ according to the types of vertices in the phantom edges of $M$. We describe two partitions of $M$, and two partitions of $\binom{M}{2}$ (noting that some sets of the partition could be empty). In the first partition of $M$, $\mathcal{E}_{3/9}(M,x)$, we differentiate according to how vertices relate to the vertex $x$, and to $V_{3/9}$. More specifically, for each phantom edge in $M$, we wish to distinguish whether each vertex in the phantom edge is either in $V_{4/9} \cup V_{5/9}$ or not, and if not, then, relative to $x$, whether each vertex in the phantom edge has smaller or larger index than $x$:
$$\begin{array}{rcl}
    E_{\ell \ell \ell}(M, x) & := & \{e \in M|\ e_1, e_2, e_3\  \mbox{are} \ x\mbox{-little}\}\\ 
    E_{\ell \ell L}(M, x) & := & \{e \in M|\ e_1 \ \mbox{and}\  e_2\ \mbox{are}\  x\mbox{-little}; \ e_3\ \mbox{is}\  x\mbox{-large}\}  \\
    E_{\ell \ell N}(M,x) & := & \{e \in M|\ e_1 \ \mbox{and}\  e_2\ \mbox{are}\ x\mbox{-little}; e_3\ \mbox{is not-small}\} \\ 
    E_{\ell L L}(M, x) & := & \{e \in M|\ e_1 \ \mbox{is}\  x\mbox{-little};\  e_2 \ \mbox{and}\  e_3\  \mbox{are}\  x\mbox{-large}\} \\
    E_{\ell L N}(M,x) & := & \{e \in M|\ e_1 \ \mbox{is}\  x\mbox{-little};\  e_2 \ \mbox{is}\  x\mbox{-large};\  e_3\  \mbox{is not-small}\} \\
    E_{\ell N N}(M,x) & := & \{e \in M|\ e_1 \ \mbox{is}\  x\mbox{-little}; \ e_2 \ \mbox{and}\  e_3 \ \mbox{are not-small}\} \\
    E_{L L L}(M, x) & := & \{e \in M|\ e_1, e_2, e_3\  \mbox{are}\  x\mbox{-large}\} \\
    E_{LLN}(M,x) & := & \{e \in M|\ e_1 \ \mbox{and}\  e_2 \ \mbox{are}\  x\mbox{-large};\  e_3\  \mbox{is not-small}\} \\
    E_{LNN}(M,x) & := & \{e \in M|\ e_1 \ \mbox{is}\  x\mbox{-large};\  e_2 \ \mbox{and}\  e_3 \ \mbox{are not-small}\} \\
    E_{N N N}(M,x) & := & \{e \in M|\ e_1, e_2, e_3\  \mbox{are not-small}\} \\
\end{array}$$
When it is clear from context, we drop the $(M,x)$. Note that we use $\ell$ to denote $x$-little vertices, $L$ to denote $x$-large vertices, and $N$ to denote not-small vertices, i.e. those in $V_{4/9} \cup V_{5/9}$. We say a vertex $v$ is {\it of type $\ell$ (with respect to $x$)} if $v$ is $x$-little. Similarly, we say that $v$ is {\it of type $L$ (with respect to $x$)} if $v$ is $x$-large, and that $v$ is {\it of type $N$} if $v$ is not-small. For $F \in \mathcal{E}_{3/9}(M,x)$ and $e \in F$, we say that $e$ is \emph{of type $F$}. We take the convention that we order the vertices in a phantom edge according to the following total order on their vertex type with respect to $x$: \begin{equation}\label{eq:order}
    \ell ~ < ~L~ <~ ~N.
\end{equation}
We also take the natural partial ordering on the sets of the partition $\mathcal{E}_{3/9}(M,x)$ acquired from \eqref{eq:order}, that is, we take the product of the linear orders. For example, given $M$ and $x$, we have that 
$$E_{\ell \ell \ell} < E_{\ell\ell N} < E_{\ell L N} < E_{N N N},$$ 
but $E_{\ell N N}$ and $E_{L L L}$ are incomparable. We also extend this ordering to the elements of the sets in $\mathcal{E}_{3/9}(M,x)$. That is, for $e \in E_{\ell \ell \ell}$, $e' \in E_{\ell\ell N}$, $e'' \in E_{\ell L N}$ and $e''' \in E_{N N N}$, we have $e<e'<e''<e'''$, and for $f \in  E_{\ell N N}$ and $f' \in E_{L L L}$,
we have both that $f \nless f'$ and $f \ngtr f'$. We shall use this partition to understand when we may obtain large and $4/9$th swapping pairs.

In the second partition of $M$, denoted by $\mathcal{E}_{5/9}(M)$, we differentiate based on the number of vertices from $V_{5/9}$. This is used in order to understand when we may obtain $5/9$th swapping pairs:
$$\begin{array}{rcl}
    E_{bbb}(M) & := & \{e \in M|\ e_1, e_2, e_3\  \mbox{are not big}\} \\
    E_{bbB}(M) & := & \{e \in M|\ e_1, e_2\  \mbox{are not big};\  e_3\  \mbox{is big}\} \\
    E_{bBB}(M) & := & \{e \in M|\ e_1 \ \mbox{is not big};\  e_2, e_3 \ \mbox{are big}\}\\ 
    E_{BBB}(M) & := & \{e \in M|\ e_1, e_2, e_3\  \mbox{are big}\} 
\end{array}$$
Note that we use $B$ to denote $5/9$th (big) vertices, and $b$ to denote all vertices which are not in $V_{5/9}$. We say that a vertex $v$ is {\it of type $b$} if $v$ is not big, {\it of type $B$} if $v$ is big, and we have the total ordering $b < B$ on these two vertex types. This ordering extends to a total ordering on the sets in $\mathcal{E}_{5/9}(M)$ given by 
$$E_{bbb}(M) < E_{bbB}(M) < E_{bBB}(M) < E_{BBB}(M).$$
As before, we also consider the related ordering to the elements of the sets, so that we may write $e<e'<e''<e'''$ when $e \in E_{bbb}(M)$, $e' \in E_{bbB}(M)$, $e'' \in E_{bBB}(M)$ and $e''' \in E_{BBB}(M)$.

Let $\{E_1, E_2\} \in \mathcal{E}_{3/9}(M,x)^{(2)}$, and $\{F_1, F_2\} \in \mathcal{E}_{5/9}(M)^{(2)}$ (where repetition is allowed). We call a pair of phantom edges $\{e,f\} \in \binom{M}{2}$ \emph{type $E_1E_2$ for $x$} if there exists $i, j \in \{1,2\}$ with $i \neq j$ such that $e \in E_i$ and $f \in E_j$ and \emph{type $F_1F_2$} if $e \in F_i$ and $f \in F_j$.

We allow the orders on $\mathcal{E}_{3/9}(M,x)$ and $\mathcal{E}_{5/9}(M)$ to extend in the natural way to $\mathcal{E}_{3/9}(M,x)^{(2)}$ and $\mathcal{E}_{5/9}(M)^{(2)}$ respectively\footnote{That is, for $E_1, E_2, F_1, F_2 \in \mathcal{E}_{3/9}(M,x)$ ($\mathcal{E}_{5/9}(M)$), we have that $E_1E_2 \leq F_1F_2$ if and only if there exists $i,j \in \{1,2\}$ and $k, \ell \in \{1,2\}$ with $i \neq j$ and $k \neq \ell$ such that $E_i \leq F_k$ and $E_j \leq F_{\ell}$.}, and define two partitions of $\binom{M}{2}$ according to these partitions. The first partition is required for Lemma~\ref{lem:3/9ths}, while the second is required for Lemma~\ref{lem:4/9ths} (both stated later). \\~\\ 
\textbf{Partition 1:}
$$\begin{array}{lll}
    T^4_{M,x} & := & \{\{e,f\} \in \binom{M}{2}|\  \{e,f\} \ \mbox{is of type}\  E_1E_2 \geq F_1F_2 \mbox{~for some pair~} \\
     & & F_1F_2 \in \{E_{\ell \ell \ell}E_{NNN}, E_{LLL}E_{LLL}\}  \ \mbox{for}\  x \}\\\\
    T^5_{M,x} & := & \{\{e,f\} \in \binom{M}{2} \setminus T^4_{M,x} | \ \{e,f\} \ \mbox{is of type}\  E_1E_2 \geq F_1F_2 \mbox{~for some pair~} \\
    & & F_1F_2 \in \{E_{\ell\ell N}E_{\ell N N}, E_{\ell \ell L}E_{L N N}, E_{\ell L N}E_{\ell L N}, E_{\ell LL}E_{LLN}\}   \ \mbox{for}\  x\}\\\\
    T^6_{M,x} & := & \{\{e,f\} \in \binom{M}{2} \setminus (\bigcup_{i =4}^{5} T^i_{M,x}) | \ \{e,f\} \ \mbox{is of type}\  E_1E_2 \geq F_1F_2 \mbox{~for some} \\  
    & & \mbox{pair~} F_1F_2 \in \{E_{\ell \ell N}E_{\ell \ell N}, E_{\ell \ell L}E_{LLN}\}   \ \mbox{for}\  x\} \\\\
    T^7_{M,x} & := & \{\{e,f\} \in \binom{M}{2} \setminus (\bigcup_{i =4}^{6} T^i_{M,x})| \ \{e,f\} \ \mbox{is of type}\ E_1E_2 \geq F_1F_2 \mbox{~for some} \\ 
    & & \mbox{pair~} F_1F_2 \in \{E_{\ell \ell \ell}E_{\ell\ell N}, E_{L L L}E_{\ell\ell L}, E_{\ell LL}E_{\ell LL}\}   \ \mbox{for}\  x\} \\\\
    T^8_{M,x} & := & \{\{e,f\} \in \binom{M}{2} \setminus (\bigcup_{i =4}^{7} T^i_{M,x})| \ \{e,f\} \ \mbox{is of type}\ E_1E_2 \geq E_{\ell\ell L}E_{\ell\ell L}  \ \mbox{for}\  x\} \\
    \\
    T^{10}_{M,x} & := & \binom{M}{2} \setminus (\bigcup_{i =4}^{8} T^i_{M,x})\\
    \end{array}$$ 
\textbf{Partition 2:}
$$\begin{array}{lll}
    S_M^4 & := & \{\{e,f\} \in \binom{M}{2} | \ \{e,f\} \ \mbox{is of type}\   E_1E_2 \geq E_{bbb}E_{BBB}   \ \mbox{for}\  x\} \\
    \\
    S_M^5 & := & \{\{e,f\} \in \binom{M}{2} \setminus (S^4_{M})| \ \{e,f\} \ \mbox{is of type}\   E_1E_2 \geq E_{bbB}E_{bBB}   \ \mbox{for}\  x\} \\
    \\
    S_M^6 & := & \{\{e,f\} \in \binom{M}{2} \setminus (\bigcup_{i =4}^{5} S^i_{M})| \ \{e,f\} \ \mbox{is of type}\   E_1E_2 \geq E_{bbB}E_{bbB}   \ \mbox{for}\  x\} \\
    \\
    S_M^7 & := & \{\{e,f\} \in \binom{M}{2} \setminus (\bigcup_{i =4}^{6} S^i_{M})| \ \{e,f\} \ \mbox{is of type}\   E_1E_2 \geq E_{bbb}E_{bbB}   \ \mbox{for}\  x\} \\
    \\
    S^{10}_M & := & \binom{M}{2} \setminus (\bigcup_{i =4}^{7} S^i_{M})\\
\end{array}$$

 The motivation for the superscript $i \in \{4, 5, 6, 7, 8\}$ ($i \in \{4, 5, 6, 7\}$), is that given a phantom matching $M$ and two vertices $x \neq y \in L(M)$, we will show that $\{x,y\}$ is a swapping pair for every pair $\{e, f\} \in T^i_{M,x} \cap T^i_{M,y}$ ($\{e, f\} \in S^i_{M}$), such that $|E(L_x(e,f))|,|E(L_y(e,f))|\geq i$. Observe that $T^{10}_{M,x}$ and $T^{10}_{M,y}$ ($S^{10}_{M}$) consist(s) of the pairs $\{e,f\}$ such that even if $L_x(e,f)$ and $L_y(e,f)$ were complete bipartite graphs, $\{x,y\}$ would not be a swapping pair for $\{e,f\}$.
    
    Before moving onto our results relating to the various partitions defined above, for vertices $x \neq y \in L(M)$ we additionally define $\mathcal{E}_{3/9}(M,x,y)$ to be the collection containing the following sets
    $$\begin{array}{rcl}
    E_{\ell \ell \ell}(M, x, y) & := & \{e \in M|\ e_1, e_2, e_3\  \mbox{are} \ x\mbox{-little and} \ y\mbox{-little}\}\\ 
    E_{\ell \ell L}(M, x, y) & := & \{e \in M|\ e_1 \ \mbox{and}\  e_2\ \mbox{are}\  x\mbox{-little and} \ y\mbox{-little}; \\ & & e_3\ \mbox{is}\  x\mbox{-large and} \ y\mbox{-large}\}  \\
    E_{\ell \ell N}(M,x, y) & := & \{e \in M|\ e_1 \ \mbox{and}\  e_2\ \mbox{are}\ x\mbox{-little and} \ y\mbox{-little}; \\ & & e_3\ \mbox{is not-small}\} \\ 
    E_{\ell L L}(M, x, y) & := & \{e \in M|\ e_1 \ \mbox{is}\  x\mbox{-little and} \ y\mbox{-little}; \\ & &  e_2 \ \mbox{and}\  e_3\  \mbox{are}\  x\mbox{-large and} \ y\mbox{-large}\} \\
    E_{\ell L N}(M,x, y) & := & \{e \in M|\ e_1 \ \mbox{is}\  x\mbox{-little and} \ y\mbox{-little};\  e_2 \ \mbox{is}\  x\mbox{-large and} \ y\mbox{-large};\\ & &  e_3\  \mbox{is not-small}\} \\
    E_{\ell N N}(M,x, y) & := & \{e \in M|\ e_1 \ \mbox{is}\  x\mbox{-little and} \ y\mbox{-little}; \ e_2 \ \mbox{and}\  e_3 \ \mbox{are not-small}\} \\
    E_{L L L}(M, x, y) & := & \{e \in M|\ e_1, e_2, e_3\  \mbox{are}\  x\mbox{-large and} \ y\mbox{-large}\} \\
    E_{LLN}(M,x, y) & := & \{e \in M|\ e_1 \ \mbox{and}\  e_2 \ \mbox{are}\ x\mbox{-large and} \ y\mbox{-large};\  e_3\  \mbox{is not-small}\} \\
    E_{LNN}(M,x, y) & := & \{e \in M|\ e_1 \ \mbox{is}\  x\mbox{-large and} \ y\mbox{-large};\  e_2 \ \mbox{and}\  e_3 \ \mbox{are not-small}\} \\
    E_{N N N}(M,x,y) & := & \{e \in M|\ e_1, e_2, e_3\  \mbox{are not-small}\} \\
\end{array}$$
    
    Note that this is not necessarily a partition of $M$. We introduce this definition to avoid any ambiguity later. Given $E_1E_2 \in \mathcal{E}_{3/9}(M,x,y)^{(2)}$, we say that $\{e, f\} \in \binom{M}{2}$ is {\it of type $E_1E_2$ for $x$ and $y$} if there exist $i, j \in \{1,2\}$ with $i \neq j$ such that $e \in E_i$ and $f \in E_j$. 
    
\subsection{Results}

\begin{prop}\label{proposition:xymatching}
Let $H$ be a 3-graph and $e,f$ be disjoint $3$-sets of vertices in $V(H)$. Let $x,y \in V(H)\setminus (e \cup f)$ and $|E(L_x(e,f))|, |E(L_y(e,f))| \geq 4$. Then $\{e,f\}$ has an $\{x,y\}$-matching.
\end{prop}
\begin{proof}
Let $E^*=E(L_x(e,f)) \cap E(L_y(e,f))$. We separate into cases based on $|E^*|$. 

If $|E^*|\geq 4$, then it is straightforward to see that we can find a matching of size two on the edges of $E^*$, which yields an $\{x,y\}$-matching. 

If $|E^*|=3$ and contains two disjoint edges then we obtain an $\{x,y\}$-matching. Else, the three shared edges all share a vertex. Without loss of generality, let this vertex be $e_1 \in e$. Consider the edge in $E(L_x(e,f))\setminus E^*$. This edge does not contain $e_1$ and contains only one of the three neighbours of $e_1$ in $f$. Thus taking this edge for $x$, at least one of the edges containing $e_1$ is disjoint and can be taken for $y$, yielding an $\{x,y\}$-matching for $\{e,f\}$.

Suppose now that $|E^*|=2$. Either this is an $\{x,y\}$-matching, or the two edges share a vertex. Without loss of generality, suppose this vertex is $e_1 \in e$. Then between the four edges which are not shared there exists an edge which does not contain $e_1$. Then we find an $\{x,y\}$-matching by taking this edge, and any edge in $E^*$ which does not intersect the chosen edge. 

If $|E^*|=1$, then there are a total of six edges which are not shared. At least one of these is disjoint from the edge in $E^*$, allowing us to find an $\{x,y\}$-matching.

Finally suppose that $E^*=\emptyset$. Since $|E(L_x(e,f))|\geq 4$, it follows that $L_x(e,f)$ must have a vertex of degree at least two in $e$. Without loss of generality, let this vertex be $e_1$. Then there exists an edge in $E(L_y(e,f))$ not containing $e_1$ (as $|E(L_y(e,f))| \geq 4$). Without loss of generality let this edge be $e_2f_1$. Since $e_1$ has degree at least two in $E(L_x(e,f))$, there exists $j \in \{2,3\}$ such that $e_1f_j \in E(L_x(e,f))$. Then $\{xe_1f_j, ye_2f_1\}$ is an $\{x,y\}$-matching for $\{e,f\}$.
\end{proof}

To prove the next lemma, we use the following definitions: for a phantom matching $M$ in $H$ and phantom edges $e,f \in M$, let $e = \{e_1, e_2, e_3\}$ and $f = \{f_1, f_2, f_3\}$. Let $x \in V(H)\setminus (e \cup f)$. 
For $g \in (e \cup f)$, we say that $L_x(e,f)$ has a \emph{star at $g$} if $d_{L_x(e,f)}(g)=3$.
If there exist $i,j \in \{1,2,3\}$ such that there are stars at $e_i$ and $f_j$ in $L_x(e,f)$ then we say that $L_x(e,f)$ has a \emph{fan at $e_if_j$}.
In addition, we remind the reader that we take the convention of ordering the vertices in a phantom edge according to the order on their vertex type with respect to relevant parameters ($x$ and $y$), and splitting ties arbitrarily.

\begin{lem}\label{lem:3/9thsswapping}
Let $\gamma > 0$ and $H$ be a 3-graph on $n$ vertices as given in Theorem~\ref{thm:almostmain}. Let $M$ be a phantom matching in $H$, let $x, y \in L(M)$ be $3/9$th vertices with $I_H(x)>I_H(y)$, and consider $\{e,f\} \in \binom{M}{2}$. Suppose there exists $i \in \{4,5,6,7,8\}$ such that $\{e,f\} \in T^i_{M,x}$ and $\{e,f\} \in T^i_{M,y}$. Suppose further that there exist types $E_1, E_2$ with $\{E_1, E_2\} \in \mathcal{E}_{3/9}(M,x,y)^{(2)}$ such that $\{e,f\}$ is of type $E_1E_2$ for $x$ and $y$; and $|E(L_x(e,f))|, |E(L_y(e,f))| \geq i$. Then $\{x,y\}$ is a swapping pair for $\{e,f\}$.
\end{lem}
\begin{proof}
First note that, given $I_H(x)>I_H(y)$ and $\{e, f\}$ being of type $E_1E_2$ for both $x$ and $y$, we have that $z \in (e \cup f)$ is $x$-little ($x$-large) if and only if it is $y$-little ($y$-large). Indeed, assume for a contradiction that $z$ is both $x$-little and $y$-large. Then, in order for $\{e, f\}$ to be of type $E_1E_2$ for both $x$ and $y$, there exists $z' \in (e \cup f)$ such that $z'$ is both $x$-large and $y$-little. But $I_H(x) > I_H(y)$ and so no vertex can be both $x$-large and $y$-little. Hence no such $z$ exists. Thus throughout what follows whenever we say a vertex is $x$-little ($x$-large) we implicitly mean it is also $y$-little ($y$-large), and vice versa.

Recall that we say an $\{x,y\}$-matching is \emph{good} for $\{e,f\}$ if it is witness for $\{x,y\}$ being a swapping pair for $\{e,f\}$. Note that $\{e,f\}$ having $\{x,y\}$ as a swapping pair is monotonous with respect to the partial order, that is, if $\{e,f\}$ was instead of type $F_1F_2$ and $E_1E_2 \leq F_1F_2$, then $\{x,y\}$ would still be a swapping pair for $\{e,f\}$. We prove the lemma by considering the different cases.

{\bf Case 1: $i=4$.}
Since  $|E(L_x(e,f))|, |E(L_y(e,f))| \geq 4$, by Proposition~\ref{proposition:xymatching} we have that $\{e,f\}$ has an $\{x,y\}$-matching. 
Observe that, since $\{e,f\} \in T^4_{M,x}$ and $\{e,f\} \in T^4_{M,y}$, whichever vertices of $e$ and $f$ are in this $\{x,y\}$-matching for $\{e,f\}$, we either have that at least one of the remaining vertices is not-small, or both the remaining vertices are $x$-large and $y$-large respectively. Hence $\{x,y\}$ is a swapping pair for $\{e,f\}$.

{\bf Case 2: $i=5$.}
It suffices to prove that Lemma~\ref{lem:3/9thsswapping} holds for $\{e,f\}$ of type $E_{\ell\ell N}E_{\ell N N}, E_{\ell \ell L}E_{L N N}$, $E_{\ell L N}E_{\ell L N}$ and $E_{\ell LL}E_{LLN}$. We consider these cases one by one, recalling our convention that $e = \{e_1, e_2, e_3\}$ and $f = \{f_1, f_2, f_3\}$.

{\bf Case 2.1: $\{e,f\}$ is of type $E_{\ell\ell N}E_{\ell N N}$ for $x$ and $y$.}
By convention, $e_1, e_2, f_1$ are $x$-little and $e_3, f_2, f_3$ are not-small. If either $e_1f_1$ or $e_2f_1$ is an edge in either $L_x(e,f)$ or $L_y(e,f)$, without loss of generality assume that $e_1f_1$ is an edge in $L_x(e,f)$. Then either $L_y(e,f)$ contains an edge disjoint from $e_1f_1$, which in turn yields that $\{x,y\}$ is a $4/9$th swapping pair for $\{e,f\}$, or $L_x(e,f)=L_y(e,f)$ and they are both precisely the fan at $e_1f_1$. Then choosing an $\{x,y\}$-matching that covers $e_1$ and $f_1$, we see that $\{x,y\}$ is a $4/9$th swapping pair for $\{e,f\}$.

Else we have that neither $e_1f_1$ or $e_2f_1$ are edges in $L_x(e,f)$ and $L_y(e,f)$. If $e_3f_1$ is an edge in either $L_x(e,f)$ and $L_y(e,f)$ we see that $\{x,y\}$ is a $4/9$th swapping pair for $\{e,f\}$, since there is certainly an $\{x,y\}$-matching containing $f_1$. Otherwise $L_x(e,f)$ and $L_y(e,f)$ only have edges incident to $f_2$ and $f_3$. Since $|E(L_x(e,f))|, |E(L_y(e,f))| \geq 5$ we have that there exist $i,j \in \{1,2\}, i \neq j,$ and $k, \ell \in \{2,3\}, k \neq \ell,$ such that $\{xe_{i}f_{k}, ye_{j}f_{\ell}\}$ is a good $\{x,y\}$-matching for $\{e,f\}$.
 
{\bf Case 2.2: $\{e,f\}$ is of type $E_{\ell \ell L}E_{L N N}$ for $x$ and $y$.}
 By convention, $e_1, e_2$ are $x$-little, $e_3, f_1$ are $x$-large and $f_2, f_3$ are not-small. We showed in Case 2.1 that either we obtain a matching leaving one of $f_2$ or $f_3$ unmatched, or we obtain a matching leaving $e_3$ and $f_1$ unmatched. In this setting, the same argument yields that $\{x,y\}$ is a $4/9$th swapping pair for $\{e,f\}$ in the first case, and a large swapping pair for $\{e,f\}$ in the second case.

{\bf Case 2.3: $\{e,f\}$ is of type $E_{\ell L N}E_{\ell L N}$ for $x$ and $y$.}
 By convention, $e_1, f_1$ are $x$-little, $e_2, f_2$ are $x$-large and $e_3, f_3$ are not-small. We start as in Case 2.1. If $e_1f_1$ is an edge in either $L_x(e,f)$ or $L_y(e,f)$, without loss of generality assume that $e_1f_1$ is an edge in $L_x(e,f)$. Then either $L_y(e,f)$ contains an edge disjoint from $e_1f_1$, which in turn either swaps out at least one $4/9$th vertex, or two $x$-large vertices, yielding that $\{x,y\}$ is a swapping pair for $\{e,f\}$, or $L_y(e,f)=L_x(e,f)$ and they are both precisely the fan at $e_1f_1$. Then choosing an $\{x,y\}$-matching that covers $e_1$ and $f_1$, we see that $\{x,y\}$ is again a swapping pair for $\{e,f\}$.

So suppose neither $L_x(e,f)$ nor $L_y(e,f)$ contains the edge $e_1f_1$. Suppose that $e_1$ is not isolated in $L_x(e,f)$ and $f_1$ is not isolated in $L_y(e,f)$ (or vice versa). Then again, choosing an edge which covers $e_1$ and another to cover $f_1$ yields a good swapping pair. So in the remaining case we have that precisely one of $e_1$ or $f_1$ is isolated in both $L_x(e,f)$ and $L_y(e,f)$ (since $|E(L_x(e,f))|, |E(L_y(e,f))| \geq 5$, so we cannot have that both are isolated). Without loss of generality, assume that $e_1$ is isolated. If $L_x(e,f)=L_y(e,f)$, then we can either take $e_2f_1$ and $e_3f_2$ or $e_2f_2$ and $e_3f_1$ to yield a good swapping pair. If $L_x(e,f) \neq L_y(e,f)$ then their union is a $K_{2,3}$, and in each graph at most one of the edges of $K_{2,3}$ is missing. So $e_2f_1$ is an edge in one of the link graphs; without loss of generality let $e_2f_1 \in L_x(e,f)$. If $e_3f_2 \notin L_y(e,f)$ then we must have $e_3f_2 \in L_x(e,f)$ and also $e_2f_1 \in L_y(e,f)$. In either case we yield that $\{x,y\}$ is a swapping pair for $\{e,f\}$.

{\bf Case 2.4: $\{e,f\}$ is of type $E_{\ell LL}E_{LLN}$ for $x$ and $y$.}

 By convention, $e_1$ is $x$-little, $e_2, e_3, f_1, f_2$ are $x$-large and $f_3$ is not-small. First note that if we can cover $e_1$ by an $\{x,y\}$-matching, then we are done. We consider the number of edges in both $L_x(e,f)$ and $L_y(e,f)$ containing $e_1$. Assume $e_1$ is incident to at least 2 edges in, without loss of generality, $L_x(e,f)$. Then there exist $i,j \in \{1,2,3\}$ with $i \neq j$, and $k \in \{2,3\}$ such that $\{xe_1f_i, ye_kf_j\}$ is a good $\{x,y\}$-matching for $\{e,f\}$.

So assume there is at most one edge incident to $e_1$ in $L_x(e,f)$ and at most one edge incident to $e_1$ in $L_y(e,f)$. If $L_x(e,f)$ $(L_y(e,f))$ has an edge incident to $e_1$ then there exists $i \in \{1,2,3\}$ such that $e_1f_i \in L_x(e,f)$ $(L_y(e,f))$. Then $L_y(e,f)$ $(L_x(e,f))$ has an edge not intersecting $e_1f_i$ (as $|E(L_y(e,f))|$ $( |E(L_x(e,f))|) \geq 5$) and $\{x,y\}$ is a swapping pair for $\{e,f\}$.

Hence assume $L_x(e,f)$ and $L_y(e,f)$ only have edges incident to $e_2$ and $e_3$. Since \newline $|E(L_x(e,f))|, |E(L_y(e,f))| \geq 5$ we have that there exist $i,j \in \{2,3\}$ with $i \neq j,$ and $k, \ell \in \{1,2\}$ with $k \neq \ell,$ such that $\{xe_if_k, ye_jf_{\ell}\}$ is a good $\{x,y\}$-matching in $\{e,f\}$.

{\bf Case 3: $i=6$.}

It suffices to prove that Lemma~\ref{lem:3/9thsswapping} holds for $\{e,f\}$ of type $E_{\ell \ell N}E_{\ell \ell N}$ and $E_{\ell \ell L}E_{LLN}$.

{\bf Case 3.1: $\{e,f\}$ is of type $E_{\ell \ell N}E_{\ell \ell N}$ for $x$ and $y$.}

By convention, $e_1, e_2, f_1, f_2$ are $x$-little and $e_3, f_3$ are not-small. Since \newline$|E(L_x(e,f))|~\geq~6$ we must have that there exist $i,j \in \{1,2\}$ such that $e_if_j \in L_x(e,f)$. Then either $L_y(e,f)$ is the union of $\{e_3, f_3\}$ and the fan at $e_if_j$, or one of the edges in $E_y:=\{e_{[2]\setminus\{i\}}f_3$, $e_3f_{[2]\setminus\{j\}},e_{[2]\setminus\{i\}}f_{[2]\setminus\{j\}}\}$ is in $E(L_y(e,f))$. In the latter case, we have that $\{xe_if_j, yab\}$ is a good $\{x,y\}$-matching for $\{e,f\}$ where $ab \in E_y$. In the former, since $|E(L_x(e,f))| \geq 6$ we observe that at least one of the following edges is in $E(L_x(e,f))$: $e_{[2]\setminus\{i\}}f_{[2]\setminus\{j\}}, e_{[2]\setminus\{i\}}f_{j}, e_{i}f_{[2]\setminus\{j\}} ,e_{[2]\setminus{i}}f_3, e_3f_{[2]\setminus{j}}$. It is then easy to find a good $\{x,y\}$-matching for $\{e,f\}$ which does not include at least one of $e_3$ or $f_3$.

{\bf Case 3.2: $\{e,f\}$ is of type $E_{\ell \ell L}E_{LLN}$ for $x$ and $y$.}

By convention, $e_1, e_2$ are $x$-little, $e_3, f_1, f_2$ are $x$-large and $f_3$ is not-small. Observe that if we can find an $\{x,y\}$-matching avoiding $e_3$ or $f_3$ then we have a good $\{x,y\}$-matching for $\{e,f\}$. Suppose at least one of $L_x(e,f)$ and $L_y(e,f)$ has three edges containing $e_3$. Without loss of generality let it be $L_y(e,f)$. Then in $L_x(e,f)$, there
are at least three edges which are incident to $e_1$ or $e_2$, and so at least one of $e_1$ and $e_2$ has an edge avoiding $f_3$. Without loss of generality let $e_2$ have such an edge. Then there exists $i \in \{1,2\}$ such that both $xe_2f_i$ and $ye_3f_{[2]\setminus\{i\}}$ are edges, yielding a good $\{x,y\}$-matching for $\{e,f\}$. Else, both $L_x(e,f)$ and $L_y(e,f)$ have at most two edges containing $e_3$. In this case they both have at least four edges avoiding $e_3$ and it is possible to find an $\{x,y\}$-matching avoiding $e_3$, and hence $\{x,y\}$ is a swapping pair for $\{e,f\}$.

{\bf Case 4: $i=7$.}

It suffices to prove that Lemma~\ref{lem:3/9thsswapping} holds for $\{e,f\}$ of type $E_{\ell \ell \ell}E_{\ell\ell N}, E_{L L L}E_{\ell\ell L}$ and $E_{\ell LL}E_{\ell LL}$.

{\bf Case 4.1: $\{e,f\}$ is of type $E_{\ell \ell \ell}E_{\ell\ell N}$ for $x$ and $y$.}

By convention, $e_1, e_2, e_3, f_1, f_2$ are $x$-little and $f_3$ is not-small. Since $|E(L_x(e,f))|,$ $|E(L_y(e,f))| \geq 7$, it follows that both have at least four edges incident to $f_1$ and $f_2$. Hence there exist $i,j \in \{1,2,3\}$ with $i \neq j,$ and $k,\ell \in \{1,2\}$ with $k \neq \ell,$ such that $\{xe_if_k, ye_jf_{\ell}\}$ is a good $\{x,y\}$-matching for $\{e,f\}$.

{\bf Case 4.2: $\{e,f\}$ is of type $E_{L L L}E_{\ell\ell L}$ for $x$ and $y$.}

The same argument as in Case 4.1 yields that $\{x,y\}$ is a large swapping pair for $\{e,f\}$.

{\bf Case 4.3: $\{e,f\}$ is of type $E_{\ell LL}E_{\ell LL}$ for $x$ and $y$.}

By convention, $e_1, f_1$ are $x$-little and $e_2, e_3, f_2, f_3$ are $x$-large. Note that since we have that both $|E(L_x(e,f))|,  |E(L_y(e,f))| \geq 7$, neither $L_x(e,f)$ nor $L_y(e,f)$ contains an isolated vertex. If $e_1f_1 \in L_x(e,f)$ $(L_y(e,f))$ then, since  $|E(L_y(e,f))|$ $(|E(L_x(e,f))|) \geq 7$, there exist $i,j \in \{2,3\}$ such that $\{xe_1f_1, ye_if_j\}$ ($\{xe_if_j, ye_1f_1\}$)  is a swapping pair for $\{e,f\}$. Otherwise it is possible to take two disjoint edges, one in $L_x(e,f)$ containing $e_1$ and one in $L_y(e,f)$ containing $f_1$, and we see that $\{x,y\}$ is a large swapping pair for $\{e,f\}$.

{\bf Case 5: $i=8$.}

It suffices to prove that Lemma~\ref{lem:3/9thsswapping} holds for $\{e,f\}$ of type $E_{\ell\ell L}E_{\ell\ell L}$.
By convention, $e_1, e_2, f_1, f_2$ are $x$-little and $e_3, f_3$ are $x$-large. Since $|E(L_x(e,f))|,$ $|E(L_y(e,f))| \geq 8$, it follows that both have at least three edges avoiding both $e_3$ and $f_3$. Hence there exist $i,j \in \{1,2\}$ with $i \neq j$ and $k,\ell \in \{1,2\}$ with $k \neq \ell$ such that $\{xe_if_k, ye_jf_{\ell}\}$ is a good $\{x,y\}$-matching for $\{e,f\}$.
\end{proof}

\begin{lem}\label{lem:fewnonlinkedges}
Let $H$ be a 3-graph on $n$ vertices as given in Theorem~\ref{thm:almostmain} and $M$ be a phantom matching in $H$. Let $x \in  L(M)$. Then $x$ is in at most $3\gamma\binom{n}{2}$ edges that are not of the form $xcd$ where $c$ belongs to a phantom edge $e$ in $M$, $d$ belongs to a phantom edge $f$ in $M$, and $e \neq f$.
\end{lem}
\begin{proof}
Recall that a phantom matching always has size $\lfloor{\frac{n-\gamma n}{3}}\rfloor$. Observe that all edges containing $x$ that are not of the form $xcd$, as given above,
either include at least one vertex in $L(M)$ that is not $x$ or are of the form $xe_1e_2$ where $e_1$ and $e_2$ both belong to the same phantom edge in $M$. There are at most $\gamma n^2$ of the former and at most $n - \gamma n$ of the latter. Hence Lemma~\ref{lem:fewnonlinkedges} holds.
\end{proof}

\begin{lem}\label{lem:5/9thswapping}
Let $H$ be a 3-graph on $n$ vertices as given in Theorem~\ref{thm:almostmain} and $M$ be a phantom matching in $H$. Let $x, y \in L(M)$ and $\{e,f\} \in \binom{M}{2}$. Suppose there exists $i \in \{4,5,6,7\}$ and $\{F_1,F_2\} \in \mathcal{E}_{5/9}(M)^{(2)}$ such that $\{e,f\} \in S^i_M$, where $\{e,f\}$ is of type $F_1F_2$ and $|E(L_x(e,f))|, |E(L_y(e,f))| \geq i$. Then $\{x,y\}$ is a $5/9$th swapping pair for $\{e,f\}$.
\end{lem}
\begin{proof}
Consider the following injective map $f:\mathcal{E}_{5/9}(M) \rightarrow \mathcal{E}_{3/9}(M,x)$ for any $x \in L(M)$, given by:
$$\begin{array}{rr}
f(E_{bbb}(M))=&E_{\ell \ell \ell}(M,x), \\ 
f(E_{bbB}(M))=&E_{\ell\ell N}(M,x), \\
f(E_{bBB}(M))=&E_{\ell NN}(M,x), \\ 
f(E_{BBB}(M))=&E_{NNN}(M,x).
\end{array}$$
Observe that $\{x,y\}$ is a $4/9$th swapping pair for $\{e,f\}$ of type $f(F_1)f(F_2)$ if and only if $\{x,y\}$ is a $5/9$th swapping pair for $\{e,f\}$ of type $F_1F_2$. 

From Case 1 in the proof of Lemma~\ref{lem:3/9thsswapping}, if $\{e,f\}$ is of type $E_{\ell \ell \ell}E_{NNN}$ with \newline 
$|E(L_x(e,f))|, |E(L_y(e,f))| \geq~4,$
then $\{x,y\}$ is a $4/9$th swapping pair for $\{e,f\}$. Hence, if $\{e,f\}$ is of type $E_{bbb}E_{BBB}$ with 
$|E(L_x(e,f))|, |E(L_y(e,f))| \geq 4,$ 
then $\{x,y\}$ is a $5/9$th swapping pair for $\{e,f\}$. 

Similarly, from Case 2.1 in the proof of Lemma~\ref{lem:3/9thsswapping}, if we have that $\{e,f\}$ is of type $E_{\ell\ell N}E_{\ell NN}$ with $|E(L_x(e,f))|, |E(L_y(e,f))| \geq 5,$ 
then $\{x,y\}$ is a $4/9$th swapping pair for $\{e,f\}$. Hence, if $\{e,f\}$ is of type $E_{bbB}E_{bBB}$ with $|E(L_x(e,f))|, |E(L_y(e,f))| \geq 5$, then $\{x,y\}$ is a $5/9$th swapping pair for $\{e,f\}$. 

From Case 3.1 in the proof of Lemma~\ref{lem:3/9thsswapping}, if we have that $\{e,f\}$ is of type $E_{\ell\ell N}E_{\ell \ell N}$ with $|E(L_x(e,f))|, |E(L_y(e,f))| \geq 6$, then $\{x,y\}$ is a $4/9$th swapping pair for $\{e,f\}$. Hence, if $\{e,f\}$ is of type $E_{bbB}E_{bbB}$ with $|E(L_x(e,f))|, |E(L_y(e,f))| \geq 6$, then $\{x,y\}$ is a $5/9$th swapping pair for $\{e,f\}$. 

Finally, from Case 4.1 in the proof of Lemma~\ref{lem:3/9thsswapping}, if we have that $\{e,f\}$ is of type $E_{\ell \ell \ell}E_{\ell \ell N}$ with $|E(L_x(e,f))|, |E(L_y(e,f))| \geq 7$, then $\{x,y\}$ is a $4/9$th swapping pair for $\{e,f\}$. Hence, if $\{e,f\}$ is of type $E_{bbB}E_{bbb}$ with $|E(L_x(e,f))|, |E(L_y(e,f))| \geq 7$, then $\{x,y\}$ is a $5/9$th swapping pair for $\{e,f\}$.

By the partial order on $\mathcal{E}_{5/9}(M)^{(2)}$, this suffices to prove the Lemma.
\end{proof}

We remind the reader that for $H$ a 3-graph on $n$ vertices with degree sequence $d_1 \leq \ldots \leq d_n$, there exists a bijection $I_H: V(H) \to [n]$ such that $I_H(x) = i$ implies that $d_H(x) := d_i$, which implies an ordering $1, \ldots, n$ of the vertices according to their position in the degree sequence.

\begin{lem}\label{lem:3/9ths}
Let $H$ be a 3-graph on $n$ vertices as given in Theorem~\ref{thm:almostmain} and $M$ be a phantom matching in $H$. Let $x \in  L(M) \cap V_{3/9}$. Then there exists $i \in \{4,5,6,7,8\}$ and $\{E_1, E_2\} \in \mathcal{E}_{3/9}(M,x)^{(2)}$ such that there are at least $\gamma \binom{n}{2}/500$ pairs $\{e,f\} \in T^i_{M,x}$ of type $E_1E_2$ for $x$ with $|E(L_x(e,f))| \geq i$.   
\end{lem}
\begin{proof}
Let $j \in [q]$ such that $I_H(x) = j$. We start by noting that, by Lemma~\ref{lem:fewnonlinkedges}, 
\begin{equation}\label{eq:smalllinkamount}
    |E(L_x(M))| \geq  \left(\frac{1}{3} + \gamma\right)\binom{n}{2} + jq.
\end{equation} 
Observe that one can place at most 
$$ \sum_{i \in \{4,5,6,7,8,10\}}(i-1)|T^i_{M,x}| $$ 
link edges into $E(L_x(M))$ such that there does not exist $i \in \{4,5,6,7,8\}$ and $\{e,f\} \in \binom{M}{2}$ with $\{e,f\} \in T^i_{M,x}$ and $|E(L_x(e,f))| \geq i$. Since there are at most $q$ small vertices in $H$, we have that $q \geq |V(M) \cap V_{3/9}|$, and by considering the different types of phantom edges in $M$ according to $\mathcal{E}_{3/9}(M,x)$, we see that 
\begin{align}\label{eq:n/6}
    q \geq & 3|E_{\ell \ell \ell}|+ 3|E_{\ell \ell L}|+3|E_{\ell LL}| + 3|E_{L L L}| + 2|E_{\ell \ell N}|\\ & + 2|E_{\ell L N}| + 2|E_{LLN}| + |E_{\ell NN}|+ |E_{LNN}|.\nonumber
\end{align}
Similarly, we know that the number of $x$-little vertices in $V(M)$ is at most $j-1$, and so from $\mathcal{E}_{3/9}(M,x)$, we also have that 
\begin{equation}\label{eq:j}
    j \geq 3|E_{\ell \ell \ell}| + 2|E_{\ell\ell L}| + 2|E_{\ell\ell N}| + |E_{\ell LL}| + |E_{\ell LN}| + |E_{\ell NN}|.
\end{equation}
One can construct a lower bound for $jq$ using \eqref{eq:n/6} and \eqref{eq:j} and, since $$\sum_{i \in \{4,5,6,7,8,10\}}|T^i_{M,x}|=\binom{|M|}{2} \leq \frac{1}{9}\binom{n}{2},$$ we have that 
$$\frac{1}{3}\binom{n}{2} \geq 3\sum_{i \in \{4,5,6,7,8,10\}}|T^i_{M,x}|.$$
We claim that 
\begin{equation}\label{eq:smalllinkamount2}
\sum_{i \in \{4,5,6,7,8,10\}}(i-1)|T^i_{M,x}| \leq \frac{1}{3}\binom{n}{2} + jq.
\end{equation}
Note that
$$\sum_{i \in \{4,5,6,7,8,10\}}(i-1)|T^i_{M,x}|=3\sum_{i \in \{4,5,6,7,8,10\}}|T^i_{M,x}| + \sum_{i \in \{4,5,6,7,8,10\}}(i-4)|T^i_{M,x}|.$$
We want to show that $$\sum_{i \in \{4,5,6,7,8,10\}}(i-4)|T^i_{M,x}| \leq jq.$$ 
Indeed, noting that \begin{equation}\label{eq:sumts}\sum_{i \in \{4,5,6,7,8,10\}}(i-4)|T^i_{M,x}|=|T^5_{M,x}|+2|T^6_{M,x}|+3|T^7_{M,x}|+4|T^8_{M,x}|+6|T^{10}_{M,x}|,\end{equation} 
and using the lower bound on $jq$ described above, we see that 
\begin{equation}\label{eq:seeappendixa}
jq - \sum_{i \in \{4,5,6,7,8,10\}}(i-4)|T^i_{M,x}| \geq 0, 
\end{equation}
as required (see Appendix~\ref{app2}).

We say that a pair $\{e,f\} \in \binom{M}{2}$ is \emph{good} for $x$ if $\{e,f\} \in T^i_{M,x}$ and $|E(L_x(e,f))| \geq i$ for some $i \in \{4,5,6,7,8\}$. It follows from \eqref{eq:smalllinkamount} and \eqref{eq:smalllinkamount2} that there are at least $\gamma\binom{n}{2}$ link edges in good pairs $\{e,f\} \in \binom{M}{2}$ for $x$. Since each good pair contains at most nine edges, this yields at least $\gamma\binom{n}{2}/9$ good pairs.  
Since each pair $\{e,f\}$ is one of 55 types for $x$, we have that there exists an $i \in \{4,5,6,7,8,10\}$ and $\{E_1E_2\} \in \mathcal{E}_{3/9}(M,x)^{(2)}$ such that there are at least $\gamma \binom{n}{2}/500$ pairs $\{e,f\} \in T^i_{M,x}$ of type $E_1E_2$ for $x$ with $|E(L_x(e,f))| \geq i$.
\end{proof}

\begin{lem}\label{lem:4/9ths}
Let $H$ be a 3-graph on $n$ vertices as given in Theorem~\ref{thm:almostmain} and $M$ be a phantom matching in $H$. Let $x \in  L(M) \cap V_{4/9}$. Then there exists $i \in \{4,5,6,7\}$ and $\{F_1,F_2\} \in \mathcal{E}_{5/9}(M)^{(2)}$ such that there are at least $\gamma \binom{n}{2}/200$ pairs $\{e,f\} \in S^i_M$ of type $F_1F_2$ with $|E(L_x(e,f))| \geq i$.
\end{lem}
\begin{proof}
The proof follows a similar outline to that of the proof of Lemma~\ref{lem:3/9ths}. By Lemma~\ref{lem:fewnonlinkedges}, 
\begin{equation}\label{eq:mediumlinkamount}
    |E(L_x(M))| \geq  \left(\frac{4}{9} + \gamma\right)\binom{n}{2},
\end{equation} 
and we have that one can place at most 
$$\sum_{i \in \{4,5,6,7,10\}}(i-1)|S^i_M|$$ link edges into $E(L_x(M))$ such that there does not exist $i \in \{4,5,6,7\}$ and $\{e,f\} \in \binom{M}{2}$ with $\{e,f\} \in S^i_M$ and $|E(L_x(e,f))| \geq i$. Since $|M| = \lfloor\frac{n-\gamma n}{3}\rfloor$, we have that $|E_{BBB}| = \lfloor\frac{n-\gamma n}{3}\rfloor - (|E_{bbb}|  + |E_{bbB}| + |E_{bBB}|)$. Furthermore, since there are in total between $|V_{3/9}| + |V_{4/9}| - \gamma n$ and $|V_{3/9}| + |V_{4/9}|$ small and medium vertices in $V(M)$, we have that $|V_{3/9}| + |V_{4/9}| - \gamma n \leq 3|E_{bbb}| + 2|E_{bbB}| + |E_{bBB}| \leq |V_{3/9}| + |V_{4/9}|$. Thus there exists 
\begin{equation}\label{eq:k}
    -\frac{\gamma n}{3} \leq k \leq \frac{2\gamma n}{3},
\end{equation}
such that $3|E_{bbb}| + 2|E_{bbB}| + |E_{bBB}| + k = |V_{3/9}| + |V_{4/9}| - \gamma n/3 \leq \frac{n-\gamma n}{3}$. Hence
\begin{equation}\label{eq:z}
    |E_{BBB}| \geq 2|E_{bbb}| + |E_{bbB}| + k.
\end{equation}
Now, since $\sum_{i \in \{4,5,6,7,10\}}|S^i_{M}|=\binom{|M|}{2} \leq \frac{1}{9}\binom{n}{2}$, we have that 
$$\frac{4}{9}\binom{n}{2} \geq 4\sum_{i \in \{4,5,6,7,10\}}|S^i_{M}|.$$ 
We claim that 
\begin{equation}\label{eq:mediumlinkamount2}
\sum_{i \in \{4,5,6,7,10\}}(i-1)|S^i_{M,x}| \leq \left(\frac{4}{9}+\frac{\gamma}{2}\right)\binom{n}{2}.
\end{equation}
To see this, note that
$$\sum_{i \in \{4,5,6,7,10\}}(i-1)|S^i_{M}|=4\sum_{i \in \{4,5,6,7,10\}}|S^i_{M}| + \sum_{i \in \{4,5,6,7,10\}}(i-5)|S^i_{M}|,$$
and thus if we can show that \begin{equation}\label{eq:gammaover2}\sum_{i \in \{4,5,6,7,10\}}(i-5)|S^i_{M}| \leq \frac{\gamma}{2}\binom{n}{2},\end{equation} then \eqref{eq:mediumlinkamount2} holds. 
Indeed, noting that \begin{equation}\label{eq:i-5}\sum_{i \in \{4,5,6,7,10\}}(i-5)|S^i_{M}|=-|S^4_{M}|+|S^6_{M}|+2|S^7_{M}|+5|S^{10}_{M}|,\end{equation}
and using \eqref{eq:k} and \eqref{eq:z}, we find that \eqref{eq:gammaover2} holds (see Appendix~\ref{app3}).

We say that a pair $\{e,f\} \in \binom{M}{2}$ is \emph{good} if $\{e,f\} \in S^i_{M}$ and $|E(L_x(e,f))| \geq i$ for some $i \in \{4,5,6,7\}$. It follows from \eqref{eq:mediumlinkamount} and \eqref{eq:mediumlinkamount2}, that there are at least $\gamma\binom{n}{2}/2$ edges in good pairs $\{e,f\} \in \binom{M}{2}$ for $x$. Since each good pair contains at most nine edges, this yields at least $\gamma\binom{n}{2}/18$ good pairs. Since each pair $\{e,f\}$ is one of 10 types, we have that there exists an $i \in \{4,5,6,7,10\}$ and $\{E_1E_2\} \in \mathcal{E}_{5/9}(M)^{(2)}$ such that there are at least $\gamma \binom{n}{2}/200$ pairs $\{e,f\} \in S^i_{M}$ of type $E_1E_2$ with $|E(L_x(e,f))| \geq i$.
\end{proof}

We are now in a position to prove Lemma~\ref{lem:betterdegree}:

\begin{proof}[Proof of Lemma~\ref{lem:betterdegree}]
Let $M$ be a phantom matching, such that $|L(M) \cap V_{5/9}|$ is as large as possible. Suppose that $|L(M) \cap V_{5/9}|=: r \leq \frac{2\gamma n}{75}$. Let $s = \frac{2\gamma n}{75}+1-r$. We wish to update $M$ to a (phantom) matching $M^*$ such that the number of matching edges does not decrease, and we have swapped out at least an additional $s$ vertices of $V_{5/9}$ so that $|L(M^*) \cap V_{5/9}|\geq s+ r > \frac{2\gamma n}{75}$.

Suppose first that $|L(M) \cap V_{4/9}| \geq \gamma n/25$. By Lemmas~\ref{lem:5/9thswapping} and \ref{lem:4/9ths}, and taking $n$ sufficiently large, since 
$$\frac{\gamma n}{25}\cdot \frac{\gamma \binom{n}{2}}{200} \geq \binom{\lfloor\frac{n-\gamma n}{3}\rfloor}{2},$$ 
there exist two $4/9$th vertices $x,y \in L(M)$ and $\{e,f\} \in \binom{M}{2}$ such that $\{x,y\}$ is a $5/9$th swapping pair for $\{e,f\}$. So we may update $M$ by swapping $e$ and $f$ for a disjoint pair of edges $e'$ and $f'$ containing $x$ and $y$ respectively, such that $((e \cup f) \setminus (e' \cup f')) \cap V_{5/9} \neq \emptyset$. Then the updated phantom matching has at least the same number of matching edges, and we have an updated leave which loses two $4/9$th vertices but gains at least one $5/9$th vertex. Hence, if there exist, say, at least $3 \gamma n/25$ vertices in $L(M) \cap V_{4/9}$, then one may apply Lemmas~\ref{lem:5/9thswapping} and \ref{lem:4/9ths} and this process at most $s \leq  \frac{2\gamma n}{75}+1$ times to get an updated phantom matching $M^*$ with leave $L(M^*)$ containing strictly more than $\frac{2\gamma n}{75}$ vertices in $V_{5/9}$, as required.

If $|L(M) \cap V_{5/9}|\leq \frac{2\gamma n}{75}$ and $|L(M) \cap V_{4/9}|\leq 3\gamma n/25$, we have that $|L(M) \cap V_{3/9}|\geq 4\gamma n/5$. By Lemmas~\ref{lem:3/9thsswapping} and \ref{lem:3/9ths}, since $|L(M) \cap V_{3/9}|\geq \gamma n/25$ and taking $n$ sufficiently large, we have that 
$$\frac{\gamma n}{25}\cdot \frac{\gamma \binom{n}{2}}{500} \geq \binom{\lfloor\frac{n-\gamma n}{3}\rfloor}{2},$$ 
and hence there exist two small vertices $x,y \in L(M)$ and $\{e,f\} \in \binom{M}{2}$, such that $\{x,y\}$ is a swapping pair for $\{e,f\}$. So we may update $M$ by swapping $e$ and $f$ for a pair of disjoint edges $e'$ and $f'$ containing $x$ and $y$ respectively, such that one of the following holds:
\begin{enumerate}
    \item $((e \cup f) \setminus (e' \cup f')) \cap V_{5/9} \neq \emptyset$,
    \item $((e \cup f) \setminus (e' \cup f')) \cap V_{4/9} \neq \emptyset$,
    \item $((e \cup f) \setminus (e' \cup f')) \cap V_{3/9} = \{e_3, f_3\}$ with $I_H(e_3) > I_H(x)$ and $I_H(f_3) > I_H(y)$.
\end{enumerate}
Then the updated phantom matching has at least the same number of matching edges, and we have an updated leave, which loses two $3/9$th vertices, and either gains at least one $5/9$th vertex, at least one $4/9$th vertex, or two new $3/9$th vertices which have strictly larger indices. As long as there are at least $\gamma n/25$ small vertices in the leave, we can do one of these swaps. Eventually we must end up with at least $3\gamma n/25$ medium vertices, or at least $\lfloor \frac{2\gamma n}{75}\rfloor +1$ big vertices. In the latter case we are done, and in the former we know that we can continue by swapping medium vertices to big vertices until we have at least $\lfloor \frac{2\gamma n}{75}\rfloor +1$ big vertices, as in the previous case.
\end{proof}

\begin{proof}[Proof of Theorem~\ref{thm:almostmain}]\label{proof:almostmainthm}
To prove this requires no further ideas than those in the proof of Lemma~\ref{lem:betterdegree}. In particular, we simply show that we can continue the swapping process until we obtain $M$ with $|L(M) \cap V_{5/9}| \geq \frac{2}{3}|L(M)|$, as required. By Corollary~\ref{cor:betterdegree} we have a matching $M$ such that $|L(M) \cap V_{5/9}|>\frac{2\gamma n}{75}$ and $|M| = \lfloor \frac{n-\gamma n}{3}\rfloor$. Now suppose that $|L(M) \cap V_{5/9}| < \frac{2}{3}|L(M)|$. From the proof of Lemma~\ref{lem:betterdegree}, if we have either at least $\gamma n/25$ small vertices or at least $\gamma n/25$ medium vertices we can continue to swap until we either gain an additional medium vertex or an additional big vertex. No swap ever reduces the number of big vertices, or the number of edges in the matching, though it may reduce the number of medium vertices in order to obtain a big vertex. Since $|L(M)|=n-3|M|\geq \gamma n$ and $|L(M) \cap V_{5/9}| < \frac{2}{3}|L(M)|$, we have $|L(M) \setminus V_{5/9}|\geq \frac{1}{3}|L(M)| \geq \frac{\gamma n}{3} - 1$. Hence we have at least $\gamma n/25$ small vertices or at least $\gamma n/25$ medium vertices. Thus we may repeatedly swap until we obtain a matching $M^*$ with $|M^*|=|M|= \lfloor \frac{n-\gamma n}{3}\rfloor$, and $|L(M^*) \cap V_{5/9}| \geq \frac{2}{3}|L(M^*)|$, as required.
\end{proof}

\section{Absorbing} \label{sec:absorbing}

In this section we prove our absorbing lemma. The proof follows very closely that of \cite[Lemma 2.4]{hps}, but our degree sequence means that we cannot use their result directly as a black box, and hence we include the proof here for completeness. 

Whilst our absorbing lemma is focused only on $3$-graphs and minimum vertex degree, we note that \cite[Lemma 2.4]{hps} gives a more general absorbing lemma for $k$-graphs defined by their minimum $t$-degree. However, the same ideas here allow for possible adaptations enabling it for use in proving further degree sequence results in other $k$-graphs. See Section~\ref{sec:comments} for further discussion on this.

\begin{lem}\label{lem:absorbing}
Let $\frac{1}{2000} \geq \gamma > 0$. There exists $n' = n'(\gamma) \in \mathbb{N}$ such that the following holds. Suppose $H$ is a 3-graph on $n \geq n'$ vertices and $q \in [(1-\sqrt{\frac{2}{3}})n]$ with degree sequence $d_1 \leq \ldots \leq d_n$ such that 
$$d_i \geq \begin{cases}
    (\frac{1}{3} + \gamma) \binom{n}{2}  & \mbox{if}\  1 \leq i \leq q,\\
    (\frac{4}{9} + \gamma)\binom{n}{2}  & \mbox{if}\  q < i \leq \frac{n}{3}, \\
   \frac{5}{9}\binom{n}{2} &  \mbox{if}\  \frac{n}{3} < i. 
    \end{cases}$$ 
Then there exists a matching $M$ in $H$ of size at most $\gamma^4n/3$ such that for any set $W \subseteq V(H)\setminus V(M)$ with $|W| \in 3\mathbb{Z}$, $|W| \leq \gamma^8n$ and $|W \cap V_{5/9}| \geq \frac{2}{3}|W|$, there exists a matching covering precisely the vertices $V(M) \cup W$.
\end{lem}

Note that here we do not require the additional $\gamma\binom{n}{2}$ in the vertices with degree at least $\frac{5}{9}\binom{n}{2}$, as in Theorem~\ref{thm:main}.\footnote{In fact, as long as two thirds of the vertices have degree at least $(\frac{5}{9}-\frac{\gamma}{5})\binom{n}{2}$, the proof below goes through without any modification.}

\begin{proof}
Let $T \in \binom{V(H)}{3}$. We say that a set $A \in \binom{V(H)}{6}$ is an absorbing set for $T$ if there exists a matching of size two in $H[A]$ and a matching of size three in $H[A \cup T]$.

\begin{prop}\label{proposition:absorbing}
For every $T \in \binom{V(H)}{3}$ with at least two vertices in $V_{5/9}:=V_{5/9}(H,0)$ there are at least $\frac{\gamma^3}{400}\binom{n}{2}^3$ absorbing sets for $T$.
\end{prop}

\begin{proof}
Let $T = \{v_1, v_2, v_3\}$ and fix this ordering of vertices in $T$. Without loss of generality let $v_2, v_3 \in V_{5/9}$. There are at most $2n$ edges which contain $v_1$ and either $v_2$ or $v_3$. Furthermore, since there are at most $q \leq (1-\sqrt{\frac{2}{3}})n$ small vertices in $H$ we have that there are at most $$q(n-q)+\binom{q}{2} \leq \left(1-\sqrt{\frac{2}{3}}\right)\sqrt{\frac{2}{3}}n^2 + \binom{(1-\sqrt{\frac{2}{3}})n}{2} \leq \frac{1}{3}\binom{n}{2}$$ 
edges containing $v_1$ and at least one small vertex. Hence, since $v_1$ has degree at least $\left(\frac{1}{3} + \gamma\right)\binom{n}{2}$ and $n$ is sufficiently large there are at least $\gamma\binom{n}{2}-2n \geq \frac{\gamma}{2}\binom{n}{2}$ edges containing $v_1$ and not containing $v_2$, $v_3$ or any small vertex. Fix such an edge $\{v_1, u_2, u_3\}$. Since $v_2 \in V_{5/9}$ and $u_2 \in V_{4/9}(H, \gamma) \cup V_{5/9}$, we have that there exist at least $\gamma \binom{n}{2} - 5n \geq \frac{3\gamma}{4} \binom{n}{2}$ pairs of vertices $\{a_2, b_2\}$ such that $a_2, b_2 \notin \{v_1,v_2,v_3,u_2,u_3\}$ and $\{v_2, a_2, b_2\}$ and $\{a_2, b_2, u_2\}$ are edges in $H$. Fix $\{a_2, b_2\}$. Then similarly there are at least $\frac{3\gamma}{4} \binom{n}{2}$ pairs of vertices $\{a_3, b_3\}$ such that $a_3, b_3 \notin \{v_1,v_2,v_3,u_2,u_3,a_2,b_2\}$ and $\{v_3, a_3, b_3\}$ and $\{a_3, b_3, u_3\}$ are edges in $H$. Thus there are at least 
$$\frac{\gamma\binom{n}{2}}{2}\cdot \left(\frac{3\gamma \binom{n}{2}}{4}\right)^2 \geq \frac{\gamma^3}{4}\binom{n}{2}^3$$ 
ordered collections which make up absorbing sets for $T$. This may be an over count of the absorbing sets themselves, but no set can be counted more than $90$ ($=\frac{6!}{2^3}$) times. Hence we obtain at least $\frac{\gamma^3}{400}\binom{n}{2}^3$ absorbing sets for $T$.
\end{proof}

Let $\mathcal{L}(T)$ denote the family of all those sets absorbing $T$. By Proposition~\ref{proposition:absorbing} we have that $|\mathcal{L}(T)| \geq \frac{\gamma^3}{400}\binom{n}{2}^3$. Choose a family of sets $\mathcal{F} \subseteq \binom{V(H)}{6}$ by selecting each set in $\binom{V(H)}{6}$ independently with probability
\begin{equation}
p=\frac{\gamma^4n}{2\binom{n}{2}^3},
\end{equation}
and note that 
\begin{equation}\label{eq:p2}
2\binom{n}{2}^3\geq 2n\binom{n}{5} \geq 12\binom{n}{6}.
\end{equation}
Then $\mathbb{E}(|\mathcal{F}|)=\binom{n}{6}p \leq \frac{\gamma^4n}{12}$, and $\mathbb{E}(|\mathcal{L}(T) \cap \mathcal{F}|)=|\mathcal{L}(T)|p \geq \frac{\gamma^7n}{800}$. Thus, by Chernoff's bound (see e.g. \cite{Chernoff}), with high probability we have that $\mathcal{F}$ has the following properties:

\begin{equation}\label{eq:f}
    |\mathcal{F}| \leq \gamma^4n/6
\end{equation}

\begin{equation}\label{eq:ltf}
    |\mathcal{L}(T) \cap \mathcal{F}| \geq \gamma^7n/1000
\end{equation} for every $T$ with at least two vertices in $V_{5/9}$.

Using \eqref{eq:p2} we have that the expected number of intersecting sets in $\mathcal{F}$ is at most $$6\binom{n}{6}\binom{n}{5} p^2 \leq \gamma^{8}n/4.$$ 
Then, by Markov's inequality, we see that with probability at least $3/4$

\begin{equation}\label{eq:intersect}
\mathcal{F} \ \mbox{contains at most} \ \gamma^8n \ \mbox{intersecting pairs of sets}.
\end{equation} 
Thus, with positive probability $\mathcal{F}$ satisfies \eqref{eq:f}, \eqref{eq:ltf} and \eqref{eq:intersect}. For each of the at most $\gamma^8n$ intersecting pairs, arbitrarily removing one from each pair, and additionally removing any set that is not an absorbing set for some $3$-set $T$, we obtain $\mathcal{F}' \subseteq \mathcal{F}$, a collection of disjoint absorbing sets such that
\begin{equation}\label{eq:ltcapf'}
    |\mathcal{L}(T) \cap \mathcal{F}'| \geq \gamma^{7}n/1000 - \gamma^8n \geq \gamma^8n
\end{equation}
for all $T$ with at least two vertices in $V_{5/9}$. Then, since $\mathcal{F}'$ is a family of pairwise disjoint absorbing sets, $H[V(\mathcal{F}')]$ contains a perfect matching, $M$, of size at most $\gamma^4n/3$. Moreover, for any $W \subseteq V(H)\setminus V(M)$ such that $|W \cap V_{5/9}| \geq \frac{2}{3}|W|$, and $\gamma^8n \geq |W| \in 3\mathbb{Z}$, we can partition $W$ into at most $\gamma^8n/3$ $3$-sets, each containing at least two vertices in $V_{5/9}$, and then absorb each one with a distinct absorbing set in $\mathcal{F}'$ by \eqref{eq:ltcapf'}. Thus $V(\mathcal{F}') \cup W$ contains a perfect matching.
\end{proof}

We now prove Theorem~\ref{thm:main}.

\begin{proof}[Proof of Theorem~\ref{thm:main}]

Without loss of generality, suppose that $0 < \gamma \leq \frac{1}{2000}$ and let $\gamma$ be fixed. Applying Lemma~\ref{lem:absorbing} yields $n'$ and letting $\gamma_1=\gamma/2$ and $\gamma_2=\gamma_1^8$, by applying Theorem~\ref{thm:almostmain} with $\gamma_1$ in place of $\gamma$, we obtain $n''$. Define $n_0 := \max\{n', \frac{n''}{1-2\gamma^4}\}$. Let $H$ be a $3$-graph on $n \geq n_0$ vertices such that $n \in 3\mathbb{Z}$ and let $q \in [(1-\sqrt{\frac{2}{3}})n]$. Suppose that $H$ has degree sequence $d_1 \leq \ldots \leq d_{n}$ satisfying 
$$d_i \geq \begin{cases}
    \left(\frac{1}{3} + \gamma\right) \binom{n}{2} + iq  & \mbox{if}\  1 \leq i \leq q,\\
    \left(\frac{4}{9} + \gamma \right)\binom{n}{2}  & \mbox{if}\  q < i \leq \frac{n}{3}, \\
   \left(\frac{5}{9} + \gamma \right)\binom{n}{2} &  \mbox{if}\  \frac{n}{3} < i. 
    \end{cases}$$
Then by Lemma~\ref{lem:absorbing}, there is a matching $M$ in $H$ of size at most $\gamma^4n/3$ such that for any set $W \subset V(H) \setminus V(M)$ with $|W| \leq \gamma^8n$, $|W| \in 3\mathbb{Z}$ and $|W \cap V_{5/9}| \geq \frac{2}{3}|W|$, we have that $H[V(M) \cup W]$ has a perfect matching. Now let $H_1:=H[V(H) \setminus V(M)]$. Then $n_1:=|V(H_1)|=n-|V(M)| \geq n(1-\gamma^4)$ and $H_1$ is a $3$-graph on $n_1$ vertices with degree sequence $d_1 \leq \ldots \leq d_{n_1}$ such that 
$$d_i \geq \begin{cases}
    \left(\frac{1}{3} + \gamma_1\right) \binom{n_1}{2} + iq  & \mbox{if}\  1 \leq i \leq q,\\
    \left(\frac{4}{9} + \gamma_1\right)\binom{n_1}{2}  & \mbox{if}\  q < i \leq \frac{n}{3}, \\
   \left(\frac{5}{9} + \gamma_1\right)\binom{n_1}{2} &  \mbox{if}\  \frac{n}{3} < i. 
    \end{cases}$$
Now, suppose that there are strictly fewer than $\frac{2n_1}{3}$ vertices $v$ satisfying 
$$d(v) \geq \left(\frac{5}{9}+ \gamma_1\right)\binom{n_1}{2}.$$
Let $r$ be the number of vertices satisfying $d(v) \geq \left(\frac{5}{9}+ \gamma_1\right)\binom{n_1}{2}$, and let $s=\frac{2n_1}{3}-r$. Note that $s \leq \gamma^4n/3$. Consider a vertex $v \in V_{4/9}(H_1, \gamma_1)$. The number of edges containing $v$ and at least two vertices in $V_{5/9}(H_1, \gamma_1)$ is at most 
$$\binom{r}{2} \leq \binom{\frac{2}{3}n_1-s}{2} \leq \binom{\frac{2}{3}n_1}{2}=\frac{4}{9}\binom{n_1}{2}-\frac{n_1}{9}.$$
Hence $v$ is in at least one edge $e$ containing at most one vertex from $V_{5/9}(H_1, \gamma_1)$. Add $e$ to $M$. We may repeat this process until we have added a set $F$ of edges to $M$ to obtain $M_1=M \cup F$, such that $H_2:= H_1[V(H_1)\setminus V(F)]$ and $|V(F)|\leq 3s \leq \gamma^4n$, and $|V_{5/9}(H_1 \setminus V(F), \gamma_1/2)| = \frac{2}{3}|V(H_2)|$. This is possible since at each iteration of the process, we either get that two-thirds of remaining vertices are in $V_{5/9}(H_1\setminus V(F), \gamma_1/2)$, or that there is a vertex in $V_{4/9}(H_1 \setminus V(F), \gamma_1/2)$ contained in an edge with at most one vertex from $V_{5/9}(H_1 \setminus V(F), \gamma_1/2)$ and hence we can repeat the process. At the end of the process we obtain $H_2$ on $n_2 \geq n_1-3s \geq n_1-\gamma^4n$ vertices with degree sequence $d_1 \leq \ldots \leq d_{n_2}$ such that 
$$d_i \geq \begin{cases}
    \left(\frac{1}{3} + 4\gamma_2\right) \binom{n_2}{2} + iq  & \mbox{if}\  1 \leq i \leq q,\\
    \left(\frac{4}{9} + 4\gamma_2\right)\binom{n_2}{2}  & \mbox{if}\  q < i \leq \frac{n_2}{3}, \\
   \left(\frac{5}{9} + \frac{\gamma}{3}\right)\binom{n}{2} &  \mbox{if}\  \frac{n_2}{3} < i, 
    \end{cases}$$
where $\left(\frac{5}{9} + \frac{\gamma}{3}\right)\binom{n}{2} \geq \left(\frac{5}{9} + 4\gamma_2\right)\binom{n_2}{2}$. Since $n_2 \geq n-2\gamma^4n \geq n_0(1-2\gamma^4) \geq n''$, and $q \leq (1-\sqrt{\frac{2}{3}})n \leq (1-\sqrt{\frac{2}{3}})(n_2+2\gamma^4n)\leq \frac{n_2}{3\sqrt{2}}$, it follows from Theorem~\ref{thm:almostmain} that $H_2$ has a matching $M_2$ covering all but at most $\gamma_2n_2+2 \leq \gamma^8n$ vertices in $H_2$. Let $L:=V(H_2) \setminus V(M_2)$. Then $|L \cap V_{5/9}(H,\gamma/3)| \geq \frac{2}{3}|L|$, and by construction we have that $L \in 3\mathbb{Z}$. By Lemma~\ref{lem:absorbing}, it follows that $H$ has a perfect matching $M_3$ on $V(M_1) \cup L$. Thus $M^*:=M_2 \cup M_3$ is a perfect matching for $H$.   
\end{proof}

\section{Concluding discussion and remarks} \label{sec:comments}

Our methods give an array of degree sequences, with a key aspect being the relationship between the number of vertices below the $4/9$ barrier, and the step-size between such vertices. In particular, our strategy requires that the step-size is at least as big as the number of vertices over which it is used, however we were unable to come up with extremal examples to demonstrate that this really is necessary. In addition, whilst the almost-perfect matching would allow up to $\frac{n}{3\sqrt{2}}$ vertices below the $4/9$ barrier, our absorbing lemma only permits at most $(1 - \sqrt{\frac{2}{3}})n$ vertices with degree below this barrier. It would be interesting to know whether this disparity can be avoided through an alternative absorbing strategy. 

\subsection{Other vertex degree sequences for 3-graphs}
In general, it would be interesting to try and push for even better degree sequences. 
Whilst our results are tight in the number of vertices at or below the $4/9$ barrier, due to evasive extremal examples, it is unclear what the maximum number of vertices that can be lowered below this is, and how far these can be lowered. 
We have given degree sequences with linear step-size, but it seems plausible that best possible sequences may instead have vertex degrees which rise in pairs and/or have a quadratic step-size. 
To illustrate what one might mean by a quadratic step-size, consider a degree sequence where we replace the step `$iq$' by something like $\sum_{j \leq i} (q-j)=iq-\frac{i(i+1)}{2}$ for some appropriate value of $q$ and $i \leq q$.\footnote{Hence why we refer to such a step-size as \emph{quadratic}.}
When looking for extremal examples, with some element of a degree sequence containing a stepping passage, some of the natural structures you might consider do have this kind of step-size.
For example, consider subsets of the vertex set $A=\{v_1,\ldots,v_{n/3}\}$ and $B=\{u_1, \ldots, u_{n/3}\}$, then adding all edges $v_iu_ju_l$ such that $i \in [n/3], j \leq i$ and $l \geq j$ yields precisely an additional neighbourhood of size $\sum_{j \leq i} (\frac{n}{3}-j)=i\frac{n}{3}-\frac{i(i+1)}{2}$ for each vertex $v_i \in A$. 
It seems feasible that these naturally arising curves in such a degree sequence might lend themselves to extremal examples as much as linear step-size would. 

In addition, it would be interesting to know whether we can find a Pos\'{a}-type degree sequence with either
\begin{enumerate}
\item[(i)] many\footnote{Where `many' denotes $cn$ vertices where $c$ is a constant \emph{not} depending on $\gamma$.} vertices below the $1/3$ barrier, or 
\item[(ii)] more than $n/3$ vertices below the $5/9$ barrier.
\end{enumerate}
For (i) it is feasible to adapt the absorbing argument, (provided we have enough vertices at the $4/9$ barrier), however the proof of Lemma~\ref{lem:3/9ths}, which enables the swapping arguments to go through, relies on vertices nearly all being asymptotically close to or above the $1/3$ barrier. For (ii), the original space barrier which shows that Theorem~\ref{thm:hps_main} is asymptotically best possible shows that you cannot have more than $2n/3$ vertices below the $5/9$ barrier.

\subsection{Exact results} Our methods rely on the additional $o(n^2)$ pairs of neighbours for each vertex in the degree sequence. Whilst the exact value of $m_1(3,n)$ is known for sufficiently large $n$, it is not entirely clear how to turn our results into exact results (i.e. omit the $\gamma \binom{n}{2}$ term). One obstacle is that the proof of the exact result for $m_1(3,n)$ does not require the additional $\gamma \binom{n}{2}$ vertex degree for the absorption part of the argument, since they use the absorbing lemma (\cite[Lemma 2.4]{hps}) which only needs minimum vertex degree at least $(1/2+\gamma)\binom{n}{2} \ll m_1(3,n)$. However, since our degree sequence includes a third of the vertices with degree below $\frac{1}{2}\binom{n}{2}$, our absorbing lemma really does require the additional $\gamma\binom{n}{2}$ degree for vertices in at least one of $V_{4/9}$ or $V_{5/9}$.

Another obstacle is that a standard method for obtaining exact results is to split into cases, namely depending on whether $H$ is close to the extremal setting or not.\footnote{This is the method used in \cite{kot_3_match}.} Since tight extremal examples evade us, it is not so clear how one might implement this idea with our current results, though our discussion for Extremal example 2 may hint at an exact degree sequence with tight extremal example.

\subsection{Degree sequence results for other hypergraphs and spanning structures}

It is natural to consider what degree sequence results may be obtained which improve on known minimum degree thresholds for perfect matchings and other spanning structures, such as Hamilton cycles and tilings of subgraphs other than $K^k_k$ in any $k$-graphs, not just for $k=3$.

Sch\"{u}lke \cite{schulke} asks about vertex degree sequences that guarantee the existence of a Hamilton cycle in a $3$-graph. Though the existence of a perfect matching does not imply the existence of a Hamilton cycle, our extremal examples containing no perfect matchings do clearly imply degree sequences for which a Hamilton cycle is not guaranteed.

Whilst for $k \geq 6$ proofs for a minimum vertex degree threshold remain elusive, there are many combinations of $t$ and $k$ for which $m_t(k,n)$ is known, and finding degree sequence results which improve on these would be extremely interesting.

\section{Acknowledgements}

The authors would like to thank Andrew Treglown for some helpful initial discussions and commenting on an earlier draft of this paper. The authors would also like to thank Stephen Gould for his early input into this project, which included coming up with the parity barrier, extremal example 2. Further, the authors are grateful to two anonymous referees for their helpful and careful reviews.

\appendix

\begin{appendices}
\section{For the proof of Lemma~\ref{lem:largermatch}}\label{app1}

In \cite[Theorem 4.4]{hps}, H\`{a}n, Person and Schacht show that assuming the largest matching in a $3$-graph $H$ with minimum degree $\delta_1(H) \geq (\frac{5}{9}+4\gamma)\binom{n}{2}$ leaves $\gamma n$ or more vertices unmatched, one can derive a contradiction (so in fact there exists a matching leaving strictly fewer than $\gamma n$ vertices unmatched). They do this by fixing a phantom matching $N$ and first noting that for any $v \in L(N)$, since $\delta_1(H) \geq (\frac{5}{9}+4\gamma)\binom{n}{2}$,
$$|E(L_v(N))| \geq \deg_H(v) - 3|N| - |L(N)|(n-|L(N)|) - \binom{|L(N)|}{2}> \left(\frac{5}{9} + \gamma\right)\binom{n}{2},$$
which is \cite[(4.1)]{hps}.
In our setting, since the minimum degree is much smaller, we cannot claim this for every vertex $v \in L(N)$, but it certainly does hold for each $v \in B(N)$ which is what we claim. H\`{a}n, Person and Schacht go on to show that either we in fact did not choose the largest matching, as we can make a switch that increases the number of matching edges in $N$, or there is a vertex $v \in L(N)$ such that the pairs $\{e,f\} \in \binom{N}{2}$ satisfying $|E(L_v(e,f))| \geq 6$ contribute at most $\gamma n^2/5$ edges to $L_v(N)$. This yields the required contradiction since, as per \cite[(4.2)]{hps}, we have 
$$|E(L_v(N))|\leq 5\binom{|N|}{2}+\gamma n^2/5 < \left(\frac{5}{9} + \gamma\right)\binom{n}{2}.$$
In our setting - the statement of Lemma \ref{lem:largermatch} - we do not assume that $M$ is a largest matching in $H$, but we show that if $M$ is not sufficiently large then we must be able to make a switch that finds a larger matching than $M$. So we may take their strategy to deduce that either there is a larger matching than $M$ (giving $M^*$ as desired), or there is a vertex $v \in L(N)$ such that the pairs $\{e,f\} \in \binom{N}{2}$ satisfying $|E(L_v(e,f))| \geq 6$ contribute at most $\gamma n^2/5$ edges to $L_v(N)$. The former is what we want to show and the latter yields a contradiction (thus yielding that the former must occur), provided that we show there is such a vertex $v \in B(N)$ (rather than just $v \in L(N)$).

The way that H\`{a}n, Person and Schacht show that there is a vertex $v \in L(N)$ such that the pairs $\{e,f\} \in \binom{N}{2}$ satisfying $|E(L_v(e,f))| \geq 6$ contribute at most $\gamma n^2/5$ edges to $L_v(N)$ is as follows. They consider four distinct subsets $Y_1, Y_2, Y_3$ and $Y_4$ of $L(N)$ and deduce that either we get a larger matching, or every vertex $v \in L(N) \setminus \bigcup_{i \in [4]} Y_i$  is a vertex such that the pairs $\{e,f\} \in \binom{N}{2}$ satisfying $|E(L_v(e,f))| \geq 6$ contribute at most $\gamma n^2/5$ edges to $L_v(N)$. Hence if $L(N) \setminus \bigcup_{i \in [4]} Y_i \neq \emptyset$ then they have found $v$ such that the contradiction occurs. Now in our case, we need that $B(N) \setminus \bigcup_{i \in [4]} Y_i \neq \emptyset$, but in fact the calculations in the proof of \cite[Theorem 4.4]{hps} show this and no extra work is required of us. In particular, Facts 4.6, 4.8, 4.10 and 4.12 show that $|Y_i| \leq \frac{\gamma n}{150}$ for each $i \in [4]$ respectively so that $\sum_{i \in [4]} |Y_i| \leq \frac{2\gamma n}{75}$. Then the details following the proof of Fact 4.12 show that any vertex $v \in L(N) \setminus \bigcup_{i \in [4]} Y_i$ satisfies $|E(L_v(N))|\leq 5\binom{|N|}{2}+\gamma n^2/5$. Now this is only a contradiction in our setting if there is a vertex in $B(N) \setminus \bigcup_{i \in [4]} Y_i$, but since we have $|B(N)|> \frac{2\gamma n}{75}$ this follows.

\section{For the proof of Lemma~\ref{lem:3/9ths}}\label{app2}

For brevity, we set 
$$\begin{array}{ccl}
    A & =& |E_{\ell, \ell, \ell}|\\
    B &= &|E_{\ell, \ell, L}|\\
    C & =& |E_{\ell, L, L}|\\
    D & =& |E_{L, L, L}|\\
    E & =& |E_{\ell, \ell, N}|\\
    F & =& |E_{\ell, L, N}|\\
    G & =& |E_{L, L, N}|\\
    H & =& |E_{\ell, N, N}|\\
    I & =& |E_{L, N, N}|\\
    J & =& |E_{N, N, N}|
\end{array}$$

By \eqref{eq:n/6} and \eqref{eq:j}, we have that 

\begin{eqnarray}\label{eq:jt}
  jq  & \geq & (3A + 2B + C + 2E + F + H ) \\ 
        & &\  \cdot \ (3A + 3B + 3C + 3D + 2E + 2F + 2G + H + I) \nonumber \\ 
        & = & 9A^2 + 15AB + 12AC + 9AD + 12AE + 9AF + 6AG + 6AH + 3AI   \nonumber \\
        & & + 6B^2 + 9BC + 6BD + 10BE + 7BF + 4BG + 5BH  + 2BI + 3C^2   \nonumber \\
        & &   + 3CD + 8CE + 5CF + 2CG + 4CH + CI + 6DE + 3DF + 3DH\nonumber \\
        & &  + 4E^2 + 6EF + 4EG + 4EH + 2EI + 2F^2 + 2FG + 3FH + FI  \nonumber \\
        & & + 2GH + H^2 + HI.\nonumber
\end{eqnarray} 

From \eqref{eq:sumts}, we have that 

\begin{eqnarray}\label{eq:sumtsfinal}
    & & \sum_{i \in \{4,5,6,7,8,10\}}(i-4)|T^i_{M,x}| \\ & = & |T^5_{M,x}|+2|T^6_{M,x}|+3|T^7_{M,x}|+4|T^8_{M,x}|+6|T^{10}_{M,x}| \nonumber\\
    & = & 6\binom{A}{2} + 6AB + 6AC + 6AD  + 3AE + 3AF + 3AG + 3AH + 3AI + 4\binom{B}{2}\nonumber \\
    & &  + 4BC  + 3BD + 3BE + 3BF  + 2BG + 3BH + BI + 3\binom{C}{2} + 3CD + 3CE \nonumber \\
    & & + 3CF + CG + 3CH + CI + 3DE + 3DF + 3DH + 2\binom{E}{2} + 2EF + 2EG \nonumber \\
    & &  + EH + EI + \binom{F}{2} + FG + FH + FI  + GH +  \binom{H}{2} + HI. \nonumber
\end{eqnarray} 

Comparing coefficients in \eqref{eq:jt} and \eqref{eq:sumtsfinal}, one can see that \eqref{eq:seeappendixa} indeed holds.

\section{For the proof of Lemma~\ref{lem:4/9ths}}\label{app3}

For brevity, we set 
$$\begin{array}{ccl}
    V & =& |E_{b, b, b}|\\
    W &= &|E_{b, b, B}|\\
    Y & =& |E_{b, B, B}|\\
    Z & =& |E_{B, B, B}|
\end{array}$$

Using \eqref{eq:i-5} and \eqref{eq:z}, we have that 

\begin{eqnarray}\label{eq:appendixbbegin}
    & & \sum_{i \in \{4,5,6,7,10\}}(i-5)|S^i_{M}| \\ 
    & = & -|S^4_{M}|+|S^6_{M}|+2|S^7_{M}|+5|S^{10}_{M}| \nonumber \\
    & = & -\binom{Z}{2} - ZY - ZW - ZV - \binom{W}{2} + 2YV + 2MV + 5\binom{V}{2} \nonumber \\
    & \leq & V + \frac{W}{2} - 3WV - W^2 - WY - \frac{3V^2}{2}  - k(3V + 2W + Y) + \frac{k}{2} - \frac{k^2}{2}.\nonumber
\end{eqnarray}

Using \eqref{eq:k} and that $3V + 2W + Y \leq \frac{n - \gamma n}{3} - k$, we have from \eqref{eq:appendixbbegin} that

\begin{eqnarray}\label{eq:appendixbend}
    & & \sum_{i \in \{4,5,6,7,10\}}(i-5)|S^i_{M}| \\ 
    & \leq & V + \frac{W}{2} - 3WV - W^2 - WY - \frac{3V^2}{2} + \frac{\gamma n^2}{9} + \frac{\gamma n}{3} \nonumber \\
    & \leq & \frac{\gamma}{2}\binom{n}{2}.\nonumber
\end{eqnarray}

\end{appendices}

\end{document}